\documentclass[a4paper,11pt]{amsart}
\usepackage{amsthm,amssymb,amsfonts,amsmath,latexsym,graphicx}

\newcommand{\beq}{\begin{eqnarray}}
\newcommand{\eeq}{\end{eqnarray}}
\newcommand{\bq}{\begin{equation}}
\newcommand{\eq}{\end{equation}}

\newcommand{\eps}{{\eps}}
\newcommand{\N}{\mathbb N}

\newcommand{\R}{{\mathbb R}}
\newcommand{\Z}{\mathbb Z}
\newcommand{\C}{\mathbb C}
\def\1{\mathbb I}
\renewcommand{\(}{\left(}
\renewcommand{\)}{\right)}
\renewcommand{\eps}{\varepsilon}
\newcommand{\W}{W}
\newtheorem{theorem}{Theorem}
\newtheorem{lemma}[theorem]{Lemma}

\newtheorem{corollary}[theorem]{Corollary}
\newtheorem{proposition}[theorem]{Proposition}

\newtheorem{remark}[theorem]{Remark}
\begin{document}

\title[Weighted Moser-Trudinger and Hardy-Sobolev inequalities]{A weighted Moser-Trudinger inequality and its relation to the Caffarelli-Kohn-Nirenberg inequalities in two space dimensions}

\author{Jean Dolbeault, Maria J. Esteban}
\address{Jean Dolbeault, Maria J. Esteban: Ceremade (UMR CNRS no. 7534), Univ. Paris-Dauphine, Pl. de Lattre de Tassigny, 75775 Paris Cedex~16, France}
\email{dolbeaul, esteban@ceremade.dauphine.fr}

\author{Gabriella Tarantello}
\address{Gabriella Tarantello: Dipartimento di Matematica. Univ. di Roma ``Tor Vergata", Via della Ricerca Scientifica, 00133 Roma, Italy}
\email{tarantel@mat.uniroma2.it}

\date{\today}

\begin{abstract}
We first prove a weighted inequality of Moser-Trudinger type depending on a parameter, in the two-dimensional Euclidean space. The inequality holds for radial functions if the parameter is larger than $-1$. Without symmetry assumption, it holds if and only if the parameter is in the interval $(-1,0]$.

The inequality gives us some insight on the symmetry breaking phenomenon for the extremal functions of the Hardy-Sobolev inequality, as established by Caffarelli-Kohn-Nirenberg, in two space dimensions. In fact, for suitable sets of parameters (asymptotically sharp) we prove symmetry or symmetry breaking by means of a blow-up method. In this way, the weighted Moser-Trudinger inequality appears as a limit case of the Hardy-Sobolev inequality.
\end{abstract}

\maketitle
\thispagestyle{empty}

%%%%%%%%%%%%%%%%%%%%%%%%%%%%%%%%%%%%%%%%%%%%%%%%%%%%%%%%%%%%%%%%%%%%%%%%%%%%%%%
%%%%%%%%%%%%%%%%%%%%%%%%%%%%%%%%%%%%%%%%%%%%%%%%%%%%%%%%%%%%%%%%%%%%%%%%%%%%%%%
\section{Introduction}\label{sect1}

By Onofri's inequality on the sphere $S^2$, see for instance \cite{Beckner-93}, we have
\bq\label{On}
\int_{S^2}e^{2\,u-2\int_{S^2}u\,d\sigma}\;d\sigma\,\leq\,e^{\Vert \nabla u \Vert^2_{L^2(S^2, d\sigma)}}\;,
\eq
for all $u\in\mathcal E=\{ u\in L^1(S^2,d\sigma)\,:\,|\nabla u|\in L^2(S^2,d\sigma)\}$, where $d\sigma$ denotes the measure induced by Lebesgue's measure on $\R^3 \supset S^2$, normalized so that $\int_{S^2}d\sigma =1$. Using the stereographic projection from $S^2$ onto $\R^2$, we see that~\eqref{On} is equivalent to the following Moser-Trudinger inequality on~$\R^2$:
\[\label{MTS2}
\int_{\R^2}e^{v-\int_{\R^2}v\,d\mu}\;d\mu\,\leq\,e^{ \frac{1}{16\,\pi}\,\Vert \nabla v \Vert^2_{L^2(\R^2, dx)}}\,,
\]
for all $v\in\mathcal D=\{ v\in L^1(\R^2,d\mu)\,:\,|\nabla v|\in L^2(\R^2,dx)\}$ where $d\mu$ denotes the probability measure
$$
d\mu=\frac{dx}{\pi\,(1+|x|^2)^2}\;.
$$
In this paper, we first generalize the above Moser-Trudinger inequality to the family of probability measures
$$
d\mu_\alpha=\frac{\alpha +1}{\pi}\,\frac{|x|^{2 \alpha}\,dx}{(1+|x|^{2\,(\alpha+1)})^2}\;,$$
for $\alpha>-1$, and investigate when the weighted inequality
\bq\label{MTalpha}
\int_{\R^2} e^{v-\int_{\R^2}v\,d\mu_\alpha}\;d\mu_\alpha\,\leq\,e^{ \frac{1}{16\,\pi\,(\alpha+1)}\,\Vert \nabla v \Vert^2_{L^2(\R^2\!,\,dx)}}\,,
\eq
holds for all $v$ in the space
$$
\mathcal E_\alpha=\Big\{v\in L^1(\R^2,d\mu_\alpha)\;:\;|\nabla v|\in L^2 (\R^2,dx)\Big\}\;.
$$
In section \ref{sect2} we prove that \eqref{MTalpha} always holds for functions in $\mathcal E_\alpha$ which are radially symmetric about the origin. Meanwhile, without symmetry assumption inequality~\eqref{MTalpha} holds in $\mathcal E_\alpha$ if and only if $\alpha\in(-1,0]$.

\medskip We use the above information to investigate possible symmetry breaking phenomena for  extremal functions of the weighted Hardy-Sobolev inequality  as established by Caffarelli-Kohn-Nirenberg (see \cite{Caffarelli-Kohn-Nirenberg-84}), in two space dimensions :
\begin{eqnarray}\label{HS}
&&\left(\int_{\R^2}\frac{|u|^p}{|x|^{bp}}\;dx\right)^{2/p}\leq\,C_{a,b}\,\int_{\R^2}\frac{|\nabla u|^2}{|x|^{2a}}\;dx\quad\forall\;u\in{\mathcal D}_{a,b}\,,\\
&&\quad\mbox{with}\;a<b\leq a+1\;,\quad p=\frac2{b-a}\;,\nonumber\\
&&\quad{\mathcal D}_{a,b}=\{|x|^{-b}\,u\in L^p(\R^2,dx)\,:\,|x|^{-a}\,|\nabla u|\in L^2(\R^2,dx)\}\;,\nonumber
\end{eqnarray}
and an optimal constant $C_{a,b}$. Typically \eqref{HS} is stated with $a<0$ (see \cite{Caffarelli-Kohn-Nirenberg-84}) so that the space ${\mathcal D}_{a,b}$ can be seen as the completion of the space $C_c^\infty(\R^2)$ of all smooth functions on $\R^2$ with compact support, with respect to the norm $\|u\|^2=\|\,|x|^{-b}\,u\,\|_p^2+\|\,|x|^{-a}\,\nabla u\,\|_2^2$. Actually \eqref{HS} holds also for $a>0$ (see section~\ref{sect2}), but in this case ${\mathcal D}_{a,b}$ is obtained as the completion with respect to $\|\cdot\|$ of the set $\{u\in C_c^\infty(\R^2)\,:\,\mbox{supp}(u)\subset\R^2\setminus\{0\}\}$. We know that for $b=a+1$, the best constant in \eqref{HS} is given by $C_{a,\,b=a+1}=a^2$ and it is never achieved (see \cite[Theorem 1.1, (ii)]{Catrina-Wang-01}). On the contrary, for $a<b<a+1$, the best constant in \eqref{HS} is always achieved, say at some function $u_{a,b}\in {\mathcal D}_{a,b}$ that we will call an {\sl extremal function,\/} but its value is not explicitly known unless we have the additional information that $u_{a,b}$ is radially symmetric about the origin. In fact, in the class of positive radially symmetric functions, the extremals of~\eqref{HS} are explicitly known (see \cite{Chou-Chu-93, Catrina-Wang-01}) and given by a multiplication by a non-zero constant and a dilation of the function
\bq\label{9.1}
u^{\rm rad}_{a,b}(x)=\Big(1+|x|^{-\frac{2a\,(1+a-b)}{b-a}}\Big)^{-\frac{b-a}{1+a-b}}\,.
\eq
See \cite{Catrina-Wang-01} for more details on existence and non-existence results and for a ``mo\-dified inversion symmetry'' property based on a generalized Kelvin transformation. Also we refer to \cite{MR2001882,Lin-Wang-04,MR2053993} for further partial symmetry results about extremal functions. On the other hand, equality is achieved by non-radially symmetric extremals for a certain range of parameters $(a,b)$ identified first in \cite{Catrina-Wang-01} and subsequently improved in \cite{Felli-Schneider-03}. In fact those results provide a rather satisfactory information about the symmetry breaking phenomenon for $u_{a,b}$ when $|a|$ is sufficiently large and also apply to any dimension $N\geq 3$, where inequality~\eqref{HS} reads as follows:
\bq\label{HSN}
\left(\int_{\R^N}\frac{|u|^p}{|x|^{bp}}\;dx\right)^{2/p}\leq\,C^N_{a,b}\int_{\R^N}\frac{|\nabla u|^2}{|x|^{2a}}\;dx\;,\quad\forall\;u\in{\mathcal D}^N_{a,b}\;,
\eq
with $p=\frac{2\,N}{(N-2)+2\,(b-a)}$\;,  ${\mathcal D}^N_{a,b}=\{|x|^{-b}\,u\in L^p(\R^N,dx)\,:\,|x|^{-a}\,|\nabla u|\in L^2(\R^N,dx)\}$, an optimal constant $C^N_{a,b}\,$, and  $a$, $b\in\R$ such that $a<(N-2)/2$, $a\leq b\leq a+1$. Again we observe that inequality~\eqref{HSN} makes sense also if $a>(N-2)/2$ and $a\leq b\leq a+1$, provided the functions are in the space ${\mathcal D}^N_{a,b}$ given by the completion with respect to $\|\cdot\|$ of the set $\{u\in C_c^\infty(\R^2)\,:\,\mbox{supp}(u)\subset\R^2\setminus\{0\}\}$.

For $N\geq 3$ and $0\leq a<(N-2)/2$, the extremal $u_{a,b}$ of \eqref{HSN} (which again exists for every $a<b<a+1$) is always radially symmetric (see \cite{Chou-Chu-93}, and for a survey on previous results see \cite{Catrina-Wang-01}). On the other hand, when $a<0$, this is ensured only in some special cases described in \cite{Lin-Wang-04,MR2053993}. Also see \cite[Theorem 4.8]{MR2001882} for an earlier but slightly less general result.

\medskip In this paper, we focus on the less investigated bidimensional case $N=2$, and besides symmetry breaking phenomena, we explore the possibility of ensuring radial symmetry (which cannot be studied as in \cite{MR2001882,Lin-Wang-04,MR2053993}) for the extremal $u_{a,b}$ according to an admissible range of  parameters $(a,b)$ (see in particular \cite[Remark 4.9]{MR2001882}).

To this purpose we check in section~\ref{Sec:HS-Extended} that \eqref{HS} (or more generally, \eqref{HSN}) holds for all $a\neq 0$ (or $a\neq(N-2)/2$ if $N\geq 3$) and not only for $a<0$ (or $a<(N-2)/2$) as it is usually found in literature. In this way we can analyze radial symmetry of the extremal $u_{a,b}$ of \eqref{HS}, in the range $a\neq 0$ and for all $b\in (a, a+1)$. We find that if $N=2$, $a\neq 0$, $b\in (a, h(a))$, with
$$
h(a)= a+\frac{|a|}{\sqrt{1+a^2}}\;,
$$
no extremal $u_{a,b}$ for \eqref{HS} is radially symmetric. This result is inspired by \cite{Felli-Schneider-03}, and it is even stated without proof for $a<0$ in \cite{Lin-Wang-04,MR2053993}. Since as $|a|\to +\infty$,
$$
0<a+1-h(a)\to 0 \;,
$$
it is reasonable to look for radially symmetric extremals when $|a|$ is small. Indeed, we will show that, if $a\to 0_+$, then $h'_+(0)=2$ (or if $a\to 0_-$, then $h'_-(0)=0$) gives the ``sharp" slope of the ratio $b/a$ that signs the transition between radial symmetry and symmetry breaking. That is, we identify two regions in the set of parameters $a$ and $b$ relative to which $u_{a,b}$ is radially symmetric, or not. The precise statement of our result is as follows (also see Figure 1 below).
%------------------------------------------------------------------------------
\begin{theorem}\label{TTT1} Let $a\neq 0$ and $N=2$.
\begin{enumerate}
\item [{\rm (i)}] If $a<b<h(a)=a+\frac{|a|}{\sqrt{1+a^2}}$, then \eqref{HS} admits only {\rm non radially symmetric} extremals.
\item [{\rm (ii)}] For every $\eps>0$, there exists $\delta>0$ such that if $|a|\in (0,\delta)$, $b\in (a,a+1)$ and either $b/a>2+\eps$ if $a>0$, or $b/a<-\eps$ if $a<0$, then the extremals of \eqref{HS} are {\rm radially symmetric,} and given by a multiplication by a non-zero constant and a dilation of the function $u^{\rm rad}_{a,b}$ defined in~\eqref{9.1}.
\end{enumerate}
\end{theorem}
%------------------------------------------------------------------------------
As a consequence of (i), we can contrast (ii) with the following statement:
\begin{enumerate}
\item [{\rm (i')}] {\it For every $\eps>0$, there exists $\delta>0$ such that if $|a|\in (0,\delta)$, $b\in (a, a+1)$ and either $b/a<2-\eps$ if $a>0$, or $b/a>\eps$ if $a<0$, then any extremal of \eqref{HS}\/} is not radially symmetric.
\end{enumerate}

\begin{figure}[ht]\includegraphics[width=7cm]{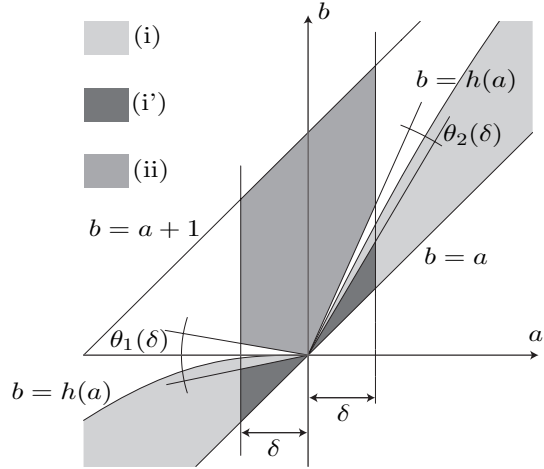}\caption{{\small\sl Radial symmetry occurs in Region {\rm (ii)}. Optimal functions are not radially symmetric in Region {\rm (i)}, and in particular in Region {\rm (i')}. The angles $\theta_1(\delta)$ and $\theta_2(\delta)$ are such that $\lim_{\delta\to 0_+}\theta_k(\delta)=0$ for $k=1$, $2$.}}\end{figure}

We will first prove (i') as a consequence of the weighted Moser-Trudinger inequality~\eqref{MTalpha}. We emphasize that such an approach makes no use of the linearized problem around the radial solution \eqref{9.1} and could be helpful in other contexts. To prove the more complete result stated in (i), we use the Emden-Fowler transformation in order to formulate~\eqref{HS} (or more generally \eqref{HSN}) as the Sobolev inequality on the cylinder $\R\times S^1$ (or more generally $\R\times S^{N-1}$). In this way we can analyze the linearized elliptic problem around the solution corresponding to \eqref{9.1} and see in which case it yields to a ``local" minimizer. We shall obtain precise informations about the linearized problem in section \ref{sect3}. This will lead us directly to the proof of (i), and will be used also to handle part (ii) of Theorem \ref{TTT1} via a blow-up analysis.

In concluding we wish to bring the reader's attention to a weighted Moser-Trudinger inequality on the cylinder $\R\times S^1$ (see Proposition~\ref{MTalphaCylinder} in section~\ref{sect5}). We believe that it helps to illustrate the nature of the symmetry breaking phenomenon analyzed here.

%%%%%%%%%%%%%%%%%%%%%%%%%%%%%%%%%%%%%%%%%%%%%%%%%%%%%%%%%%%%%%%%%%%%%%%%%%%%%%%
%%%%%%%%%%%%%%%%%%%%%%%%%%%%%%%%%%%%%%%%%%%%%%%%%%%%%%%%%%%%%%%%%%%%%%%%%%%%%%%
\section{A weighted Moser-Trudinger inequality and its connection with the weighted Hardy-Sobolev inequality}\label{sect2}

Consider the measure $\mu_\alpha$ and the Banach space $\mathcal E_\alpha$, $\alpha >-1$, defined in section~\ref{sect1}. Here and from now on, $\| v\|_2$ means $\| v\|_{L^2(\R^2,dx)}$.

\newpage
%%%%%%%%%%%%%%%%%%%%%%%%%%%%%%%%%%%%%%%%%%%%%%%%%%%%%%%%%%%%%%%%%%%%%%%%%%%%%%
\subsection{A weighted Moser-Trudinger inequality on $\R^2$}\label{sect2.1}
%------------------------------------------------------------------------------
\begin{proposition}\label{Prop:IneqCoompl} Let $\alpha>-1$. For all $v\in\mathcal E_\alpha$, there holds
\bq\label{Ineq:MTR2}
\int_{\R^2} e^{v-\int_{\R^2}v\,d\mu_\alpha}\;d\mu_\alpha\,\leq\,e^{ \frac{1}{16\,\pi\,(\alpha+1)}\,\left(\|\nabla v\|^2_2 +\alpha\,(\alpha+2)\,\|\,\frac 1{r}\,\partial_\theta v\,\|^2_2\right)}\,.
\eq
\end{proposition}
%------------------------------------------------------------------------------
\proof We use polar coordinates in $\R^2\approx\C$. For $x\in\R^2$, we let $x=r\,e^{i\theta}$, $r\geq 0$, $\theta\in [0,2\pi)$. We also consider cylindrical coordinates in $\R^3$, so that for $(y,z)\in \R^2\times\R$, we let $y=\rho\,e^{i\theta}$, $\rho\geq 0$, $\theta\in [0,2\pi)$ and $z\in\R$. In this way, we can write $\R^3\supset S^2=\{(\rho\,e^{i\theta}, z)\,:\,\rho^2+z^2=1\;\mbox{and}\;\theta\in [0, 2\pi)\}$. We recall that the inverse $\Sigma_0$ of the usual stereographic projection from $S^2$ onto~$\R^2$ is defined by
$$
\Sigma_0\big(r\,e^{i\theta}\big)=(\rho\,e^{i\theta}, z)=\left(\frac{2\,r\,e^{i\theta}}{1+r^2 }\,,\;\frac{r^{2}-1}{1+r^2}\right)\,.
$$
If $u$ is defined on $S^2$, then $v=u\,\circ\,\Sigma_0$ is defined on $\R^2$ and for any continuous real function $f$ on $\R$, we have
$$
\pi\int_{S^2}f(u)\;d\sigma=\int_{\R^2}\frac{f(v)}{(1+|x|^2)^2}\;dx\quad\mbox{and}\quad 4\pi\int_{S^2}|\nabla u|^2\;d\sigma=\int_{\R^2}|\nabla v|^2\,dx
$$
whenever $f(u)$ and $|\nabla u|^2$ belong to $L^1(S^2)$.

In order to prove the proposition, we are going to use the inverse of a {\sl dilated stereographic projection\/} given for all $\alpha>-1$ by the function $\Sigma_\alpha:\R^2\to S^2$ such that
\[\label{sigma-alpha}
\Sigma_\alpha\big(r\,e^{i\theta}\big)= \left(\frac{2\,r^{\alpha+1}\,e^{i\theta}}{1+r^{2(\alpha+1)} }\,,\;\frac{r^{2(\alpha+1)}-1}{1+r^{2(\alpha+1)}}\right)\,.
\]
Note that for any $r\geq 0$, $\theta\in [0, 2\pi)$, $\Sigma_\alpha(r\,e^{i\theta})=\Sigma_0(r^{1+\alpha}\,e^{i\theta})$ and, for any $\rho\geq 0$, $\theta\in [0, 2\pi)$ and $z \in [-1, 1]$,
$$
\Sigma_\alpha^{-1}\big((\rho\,e^{i\theta},z)\big)=\big({\textstyle \frac\rho{1-z}}\big)^{1/(\alpha+1)}e^{i\theta}\;.
$$
Now, if $f$ is a continuous real function on $\R$, $f(u),\;|\nabla u|^2\in L^1(S^2)$ and $v= u\,\circ\,\Sigma_\alpha$, then an elementary computation (see the Appendix) shows that
\begin{eqnarray*}
&&\int_{S^2}f(u)\;d\sigma =\int_{\R^2}f(v)\;d\mu_\alpha\;,\\
&&4\pi\int_{S^2}|\nabla u|^2\;d\sigma=\frac1{\alpha+1}\int_{\R^2}\Big(\,|\nabla v|^2 +\alpha\,(\alpha+2)\,\Big|\,\frac 1{r}\,\partial_\theta v\,\Big|^2\,\Big)\,dx\;.
\end{eqnarray*}
The result follows from Onofri's inequality~\eqref{On}.\endproof
Notice that we will recover Onofri's inequality as a consequence of Proposition~\ref{CT1} and the symmetry result of Theorem~\ref{TTT1}, (ii). See Remark~\ref{rem1} for details.
%------------------------------------------------------------------------------
\begin{corollary} If $\alpha\in (-1, 0]$, then \eqref{MTalpha} holds true for any $v\in\mathcal E_\alpha$.\end{corollary}
%------------------------------------------------------------------------------
\proof It is an immediate consequence of Proposition~\ref{Prop:IneqCoompl} since for $\alpha\in (-1, 0]$, we have $\alpha\,(\alpha+2)\leq 0$. \endproof

This result is optimal. While \eqref{MTalpha} remains valid for all $\alpha>-1$ among radially symmetric functions (about the origin), it fails in $\mathcal E_\alpha$ for $\alpha>0$:
%------------------------------------------------------------------------------
\begin{proposition}\label{nottrue} If $\alpha>0$, then inequality~\eqref{MTalpha} fails to hold in $\mathcal E_\alpha$.\end{proposition}
%------------------------------------------------------------------------------
\proof Let us exhibit a counter-example to \eqref{MTalpha}, which is valid for all $\alpha>0$. For any $\eps\in (0,1)$, let us consider the function $v_\eps:\R^2\to \R$ defined by
$$
2\,v_\eps = \left\{\begin{array}{ll}\log\left(\frac\eps{(\eps+\pi\,|x-\bar x|^2)^2}\right) & \mbox{if}\quad |x-\bar x|\leq 1\\
\log\left(\frac\eps{(\eps+\pi)^2}\right) & \mbox{if}\quad |x-\bar x|>1\end{array}\right.$$
where $\bar x$ denotes the point $(1,0)$. For this function we can calculate the various terms of \eqref{MTalpha}.

First we compute the l.h.s., and see that
\[\label{lhsce}
\mu_\alpha(e^{2v_\eps})=\int_{\R^2}e^{2v_\eps}\;d\mu_\alpha=I_{\alpha,\eps} +A_\alpha\,\frac\eps{(\eps+\pi)^2}
\]
where
$$
I_{\alpha,\eps}=\frac 1\eps\int_{|x-\bar x|< 1}\frac 1{\(1+\pi\,|\frac{x-\bar x}{\sqrt{\eps}}|^2\)^2}\;d\mu_\alpha
$$
and $A_\alpha=\int_{|x-\bar x|>1}d\mu_\alpha$ is finite for all $\alpha>-1$. Now, by the change of variables $x=\bar x+\sqrt\eps\,y$ and dominated convergence, we find
$$
\lim_{\eps\to0}\int_{|y|<1}\frac{|\bar x+\sqrt\eps\,y|^{2\alpha}}{\(1+|\bar x+\sqrt\eps\,y|^{2 (\alpha +1)}\)^2 \(1+\pi\,|y|^2\)^2}\;dy\;=\;\frac14 \int_{\R^2} \frac{dy}{(1+\pi\,|y|^2)^2}\;,
$$
So, for the function $v_\eps$, the l.h.s. of \eqref{MTalpha} satisfies
$$
\lim_{\eps\to0}\mu_\alpha(e^{2v_\eps})=\lim_{\eps\to0}I_{\alpha,\eps}=\frac{\alpha+1}{4\pi}\;.
$$

Next we compute the r.h.s. of \eqref{MTalpha}, that is $\frac 1{4\pi\,(\alpha+1)}\,\|\nabla v_\eps\|_2^2 +2\,\mu_\alpha(v_\eps)$ and see that
$$
\|\nabla v_\eps\|_2^2=4\,\pi\,\log\left(\frac{\eps+\pi}\eps\right)-\frac{4\,\pi^2}{(\eps+\pi)}
$$
and
$$
2\,\mu_\alpha(v_\eps)=J_{\alpha,\eps} +A_\alpha\,\log\frac\eps{(\eps+\pi)^2}\;,
$$
where
$$
J_{\alpha,\eps}=\int_{|x-\bar x|<1}\log\left(\frac\eps{ (\eps+\pi\,|{x-\bar x}|^2)^2}\right)\,d\mu_\alpha\;.
$$
Using $A_\alpha=1-\int_{|x-\bar x|<1}d\mu_\alpha$, we get
\[
2\,\mu_\alpha(v_\eps)=\log\frac\eps{(\eps+\pi)^2}+B_{\alpha,\eps}\;,\quad\! B_{\alpha,\eps}=\int_{|x-\bar x|<1}\log\left(\frac{\eps+\pi}{\eps+\pi\,|{x-\bar x}|^2}\right)^2d\mu_\alpha\;,
\]
\[
\lim_{\eps\to0}B_{\alpha,\eps}=
\int_{|x-\bar x|<1}\log\left(\frac{1}{ |{x-\bar x}|^4}
\right)\,d\mu_\alpha\;.
\]
Hence
\[
\frac 1{4\pi\,(\alpha+1)}\,\|\nabla v_\eps\|_2^2 +2\,\mu_\alpha(v_\eps)=\frac\alpha{1+\alpha}\log\eps+O(1)\quad\mbox{as}\quad\eps\to0\;,
\]
and comparing with the estimate above, we violate \eqref{MTalpha} for $\eps>0$ small enough. \endproof

%%%%%%%%%%%%%%%%%%%%%%%%%%%%%%%%%%%%%%%%%%%%%%%%%%%%%%%%%%%%%%%%%%%%%%%%%%%%%%
\subsection{The weighted Hardy-Sobolev inequality}\label{Sec:HS-Extended}
The range in which inequalities \eqref{HS} and \eqref{HSN} are usually considered can be extended as follows.
%------------------------------------------------------------------------------
\begin{lemma}\label{Lemma1} If $N=2$, then inequality~\eqref{HS} holds for any $a\neq 0$ and $b$ such that $a<b\le a+1$. If $N\geq3$, then inequality~\eqref{HSN} holds for any $a\neq (N-2)/2$ and $b$ such that $a\le b\le a+1$.
\end{lemma}
%------------------------------------------------------------------------------
\proof We use Kelvin's transformation and deal with the case $N=2$. If
$u\in {\mathcal D}_{a,b}$, then $v(x)=u\left(x/|x|^2\right)$ is such that $|x|^a\,|\nabla v|\in L^2(\R^2,dx)$. Hence, for $a>0$, $b\in (a, a+1]$, define $a'=-a$, $b'=b-2a\in (-a, -a+1]$ and apply~\eqref{HS} to the pair $(a',b')$ with $p=2/(b'-a')$ to obtain
$$
\int_{\R^2} \left( \frac{|v|^p}{|x|^{b'p}}\,dx\right)^{2/p}\,\leq\;C_{a',b'}\,\int_{\R^2} \frac{| \nabla v|^2}{|x|^{2a'}}\;dx\quad\mbox{in}\;\;{\mathcal D}_{a',b'}\;.
$$
Now, we make the change of variables $y=x/|x|^2$ and get
$$
\int_{\R^2} \left( \frac{|u|^p}{|y|^{4-b'p}}\,dy\right)^{2/p}\,\leq\;C_{a',b'}\,\int_{\R^2} \frac{| \nabla u|^2}{|y|^{-2a'}}\;dy\quad\mbox{in}\;\;{\mathcal D}_{a,b}\;.
$$
Thus we arrive at the desired conclusion with $C_{a,b}=C_{a',b'}$, since
$$
4-b'p=bp\;,\quad-2a'=2a\quad\mbox{and}\quad p=2/(b'-a')= 2/(b-a)\;.
$$
Similarly in dimension $N\geq 3$, argue as above with $a=N-2-a'$, $b\,p=2N-b'\,p$ and $p=2N/(N-2-2(b'-a'))=2N/(N-2-2(b-a))$. \endproof

Surprisingly, the case $a>0$ if $N=2$, or $a>(N-2)/2$ if $N\geq 3$, has apparently never been considered. According to our argument, it requires to define with care the space ${\mathcal D}_{a,b}$. Indeed if a function $u\in C_c^\infty(\R^N)\cap{\mathcal D}_{a,b}$ for $a>(N-2)/2$, $N\geq 2$, then $u$ must satisfy $u(0)=0$. Although optimal functions for inequality~\eqref{HSN}, $a>(N-2)/2$, $N\geq 2$, have not been studied, it has been noted in \cite[Theorem 1.4]{Catrina-Wang-01} that whenever $u>0$ satisfies the corresponding Euler-Lagrange equations, then, up to a scaling, it satisfies the ``modified inversion symmetry'' property, that is, there exists $\tau>0$ such that
\[
u(x)=\left|\frac x\tau\right|^{-(N-2-2a)}\,u\(\tau^2\,\frac x{|x|^2}\)\quad\forall\;x\in\R^N\,.
\]
The transformation $u\mapsto \left|x\right|^{-(N-2-2a)}\,u(x/|x|^2)$ is sometimes called the generalized Kelvin transformation, see {\sl e.g.\/} \cite{Chou-Chu-93}. The modified inversion symmetry formula can be shown for an optimal function $u$ using the fact that $v$ given in terms of $u$ as in the proof of Lemma~\ref{Lemma1} is also an optimal function for inequality~\eqref{HSN}, with parameters~$a'$,~$b'$.

%%%%%%%%%%%%%%%%%%%%%%%%%%%%%%%%%%%%%%%%%%%%%%%%%%%%%%%%%%%%%%%%%%%%%%%%%%%%%%
\subsection{The Moser-Trudinger inequality as a limit case of the weigh\-ted Hardy-Sobolev inequality on $\R^2$}

We now relate inequalities \eqref{MTalpha} and \eqref{HS}. In this section, we will only consider the case $a<0$. The case $a>0$ follows by Lemma~\ref{Lemma1}.

\medskip For $N=2$, $\alpha>-1$, $\eps\in(0,1)$, let us make the following special choice of parameters:
\[
a=-\frac{\eps}{1-\eps}\,(\alpha+1)\;,\quad b=a+\eps\quad\mbox{and}\quad p=\frac 2\eps\;.
\]
Let $u_\eps=u^{\rm rad}_{a,b}$ be given in \eqref{9.1}, that is
\[
u_\eps(x)=\Big(1+|x|^{2(\alpha+1)}\Big)^{-\frac\eps{1-\eps}}\;.
\]
We consider the functions
\[
f_\eps=\left[\frac{u_\eps}{|x|^{a+\eps}}\right]^{2/\eps},\quad g_\eps=\left[\frac{|\nabla u_\eps|}{|x|^a}\right]^2,
\]
and the integrals
\[
\kappa_\eps=\int_{\R^2}f_\eps\;dx\quad\mbox{and}\quad\lambda_\eps=\int_{\R^2}g_\eps\;dx\;.
\]
Straightforward computations show that
\[
\kappa_\eps = \int_{\R^2}\frac{|x|^{2\alpha}}{\big(1+|x|^{2(1+\alpha)}\big)^2}\,\frac{u_\eps^2}{|x|^{2a}}\;dx=\frac\pi{\alpha+1}\int_0^\infty\frac{s^\frac\eps{1-\eps}}{(1+s)^\frac 2{1-\eps}}\;ds\;,
\]
\[
\lambda_\eps = 4a^2\int_{\R^2}\frac{|x|^{2(2\alpha+1-a)}}{\big(1+|x|^{2(1+\alpha)}\big)^{\frac2{1-\eps}}}\,dx\;.
\]
Notice that we can use Euler's Gamma function $\Gamma(x)=\int_0^\infty s^{x-1}\,e^{-s}\,ds$, and on the basis of the well known identity:
$$
2\int_0^\infty s^{2a-1}(1+s^2)^{-b}\,ds=\frac{\Gamma(a)\,\Gamma(b-a)}{\Gamma(b)}\;,
$$
deduce for $\lambda_\eps$ the following expression:
$$
\lambda_\eps = 4\pi\,|a|\,\frac{\Gamma\big(\frac{2-\eps}{1-\eps}\big)\,\Gamma\big(\frac 1{1-\eps}\big)}{\Gamma\big(\frac 2{1-\eps}\big)}\;.
$$
%------------------------------------------------------------------------------
\begin{lemma}\label{CP1} Let $\alpha_0>-1$, $v\in C_c^\infty(\R^2)$, $w_\eps=(1+\eps\,v)\,u_\eps$. With the above notations, we have
\[
\frac 1{\kappa_\eps}\int_{\R^2}\frac{|w_\eps|^p}{|x|^{bp}}\;dx=\int_{\R^2}|1+\eps\,v|^\frac 2\eps\;\frac{f_\eps\;dx}{\int_{\R^2}f_\eps\;dx}
\label{zz2}
\]
and, as $\eps\to 0$, uniformly with respect to $\alpha\geq \alpha_0$,
\[
\int_{\R^2}\kern-7pt\frac{|\nabla w_\eps|^2}{|x|^{2a}}\,dx= \lambda_\eps+\varepsilon^2\!\left[{\textstyle \frac{8(1+\alpha)^2}{(1-\varepsilon)^2}}\kern-3pt\int_{\R^2}\kern-3pt\frac{u_\varepsilon^{2/\varepsilon}\,v}{|x|^{2(a-\alpha)}}\,dx +\!\int_{\R^2}\kern-7pt|\nabla v|^2\,\frac{u_\varepsilon^{2}}{|x|^{2a}}\,dx+O(a^2\varepsilon)\right].
\]
\end{lemma}
%------------------------------------------------------------------------------
\proof By definition of $g_\eps$, we can write
\[
\int_{\R^2}\frac{|\nabla w_\eps|^2}{|x|^{2a}}\,dx=\lambda_\eps+2\,\eps\;\underbrace{\int_{\R^2}\nabla u_\eps\cdot\nabla(u_\eps\,v)\;\frac{dx}{|x|^{2a}}}_{\rm (I)}\;+\;\eps^2\;\underbrace{\int_{\R^2}|\nabla(u_\eps\,v)|^2\,\frac{dx}{|x|^{2a}}}_{\rm (II)}\;.
\]
A simple algebraic computation shows that
\bq\label{EDO}
-\nabla\cdot\(\frac{\nabla u_\eps}{|x|^{2a}}\)=\frac{4\,a^2}\eps\,u_\eps^{\frac 2\eps-1}|x|^{2(\alpha-a)}\;.
\eq
Using \eqref{EDO} and an integration by parts, we obtain
\[
{\rm (I)}=\frac{4\,a^2}\eps\int_{\R^2}|x|^{2(\alpha-a)}\,u_\eps^{2/\eps}\,v\;dx\;.
\]
As for ${\rm (II)}$, we expand $|\nabla(u_\eps\,v)|^2$ and write
\[
{\rm (II)}=\int_{\R^2}\Big[v^2\,|\nabla u_\eps|^2+u_\eps\,\nabla(v^2)\cdot \nabla u_\eps+u_\eps^2\,|\nabla v|^2\Big]\frac{dx}{|x|^{2a}}
\]
where the first two terms can be evaluated as above using \eqref{EDO} and an integration by parts. Hence,
\[
\int_{\R^2}\Big(v^2\,|\nabla u_\eps|^2+u_\eps\nabla(v^2)\cdot \nabla u_\eps \Big)\,\frac{dx}{|x|^{2a}}=\frac{4\,a^2}\eps\int_{\R^2}|x|^{2(\alpha-a)}\,u_\eps^{2/\eps}\,v^2\;dx\;.
\]
To complete the proof we just remark that the function $|x|^{2(\alpha-a)}u_\eps^{2/\eps}$ is uniformly bounded for $\alpha\geq \alpha_0>-1$. \endproof

For a given $\alpha>-1$, we now investigate the limit as $\eps\to0$. We prove that inequality~\eqref{MTalpha} is a limiting case of inequality~\eqref{HS}, whenever~\eqref{HS} admits a radially symmetric extremal for any $\eps$ small enough. In such a case, we can write \eqref{HS} as follows:
\bq\label{IneqAsympTruc}
\frac 1{\kappa_\eps}\int_{\R^2}\frac{|w|^p}{|x|^{bp}}\;dx\leq\(\frac 1{\lambda_\eps}\int_{\R^2}\frac{|\nabla w|^2}{|x|^{2a}}\,dx\)^{1/\eps}\,.
\eq
Thus, if we take $w=w_\eps=(1+\eps\,v)\,u_\eps$, then we have:
\[
\frac 1{\kappa_\eps}\!\int_{\R^2}\kern-4pt\frac{|w_\eps|^p}{|x|^{bp}}\;dx\leq \({\textstyle 1\!\!+\!\frac{\varepsilon^2}{\lambda_\eps}\!\left[ \frac{8(1+\alpha)^2}{(1-\varepsilon)^2}\kern-3pt\int_{\R^2}\kern-3pt\frac{u_\varepsilon^{2/\varepsilon}\,v}{|x|^{2(a-\alpha)}}\,dx \!+\!\int_{\R^2}\kern-3pt\frac{|\nabla v|^2\,u_\varepsilon^{2}}{|x|^{2a}}\,dx\right]}\)^{\!1/\eps}\kern-9pt+O(a^2\varepsilon^2)\,.
\]
In particular, observe that
\[
\frac{|x|^{-bp}\,f_\eps\;dx}{\int_{\R^2}f_\eps\;dx}\sim\frac{\alpha+1}\pi\,|x|^{2\alpha}\,u_\eps^{2/\eps}\,dx\sim d\mu_\alpha(x)\quad\mbox{as}\quad\eps\to 0_+\;.
\]
%------------------------------------------------------------------------------
\begin{proposition}\label{CT1} Let us fix $\alpha>-1$ and suppose that there exists a sequence $(\eps_n)_{n\in\N}$ converging to $0$ such that the radial extremal function $u_{\eps_n}$ is also extremal for~\eqref{HS} with $(a,b,p)=(a_n,b_n,p_n)$ specified a follows,
$$
p_n=\frac2{\eps_n}\;,\quad a_n=-\frac{\eps_n}{1-\eps_n}\,(\alpha+1)\;,\quad b_n=a_n+\eps_n\;.
$$
Then the weighted Moser-Trudinger inequality~\eqref{MTalpha} holds true on $\mathcal E_\alpha$. \end{proposition}
%------------------------------------------------------------------------------
\proof As $n\to\infty$, we have
\[
\lambda_{\varepsilon_n}= 4\pi\,|a_n|+o(\eps_n)\;,\quad\kappa_{\eps_n} = \frac\pi{\alpha+1}+o(1)\;.\label{zzz1}
\]
Using Lebesgue's theorem of dominated convergence repeatedly and Lem\-ma~\ref{CP1}, for any $v\in C_c^\infty(\R^2)$ and $w_{\eps_n}=(1+\eps_n\,v)\,u_{\eps_n}$, we have
\begin{eqnarray*}
&&\hspace*{-5pt}\frac 1{\kappa_{\eps_n}}\!\int_{\R^2}\!\!\frac{|w_{\eps_n}|^{p_n}}{|x|^{b_n{p_n}}}\;dx=\!\int_{\R^2}\!\!|1+{\eps_n}\,v|^\frac 2{\eps_n}\;\frac{f_{\eps_n}\;dx}{\int_{\R^2}f_{\eps_n}\;dx}\to\!\int_{\R^2}e^{2v}\,d\mu_\alpha\;,\label{zzz2}\\
&&\hspace*{-5pt}\frac 1{\lambda_{\eps_n}}\int_{\R^2}\frac{|\nabla w_{\eps_n}|^2}{|x|^{2a_n}}\;dx=1+\eps_n\(\int_{\R^2}2v\,d\mu_\alpha+\frac{1}{4\,(1+\alpha)\,\pi}\,\|\nabla v\|^2_2\)+O(\eps_n^2)\label{zzz3}
\end{eqnarray*}
as $n\to+\infty$. The proposition follows by applying inequality~\eqref{HS} with $(a,b,p)=(a_n,b_n,p_n)$. By density we can finally choose $v$ in the larger space~$\mathcal E_\alpha$. \endproof

%------------------------------------------------------------------------------
\begin{remark}\label{rem1} Incidentally let us note that if we temporarily admit the result~{\rm (ii)} of Theorem~\ref{TTT1}, then we  find a sequence of optimal functions as required by Proposition~\ref{CT1}. In particular, for $\alpha=0$, this gives a proof of the Moser-Trudin\-ger inequality on $\R^2$ as a consequence of inequality~\eqref{HS}. Using the inverse $\Sigma_0$ of the stereographic projection, this also proves Onofri's inequality~\eqref{On} on~$S^2$. \end{remark}
%------------------------------------------------------------------------------

Let us now consider another asymptotic regime in which $\alpha\to\infty$.%------------------------------------------------------------------------------
\begin{proposition}\label{CT2} If $(\eps_n)_{n\in\N}$ and $(\alpha_n)_{n\in\N}$ are two sequences of positive real numbers such that as $n\to+\infty$,
\[
\lim_{n\to +\infty}\eps_n=0\;,\quad\lim_{n\to +\infty}\alpha_n=+\infty\quad\mbox{and}\quad a_n=-\frac{\eps_n}{1-\eps_n}\,(1+\alpha_n)\mathop{\to}_{n\to +\infty} 0_-\;,
\]
then for $n$ large enough, the radially symmetric extremal $u_{\eps_n}$ cannot be a global extremal for inequality~\eqref{HS}.
\end{proposition}
%------------------------------------------------------------------------------
\proof We argue by contradiction and assume that \eqref{IneqAsympTruc} holds with respect to the given choice of parameters. By definition of $\lambda_{\eps_n}$, $\kappa_{\eps_n}$, and Lebesgue's theorem of dominated convergence, we know that
\[\label{AA1}
\lim_{n\to +\infty}\frac{\lambda_{\eps_n}}{|a_n|}= 4\pi\quad\mbox{and}\quad\lim_{n\to +\infty}(\alpha_n+1)\,\kappa_{\eps_n}=\pi\;.
\]
If $v\in C_c^\infty(\R^2)$, then by a direct computation, we find:
\begin{eqnarray*}
&&\hspace*{-12pt}(\alpha_n+1)\int_{\R^2}\frac{|u_{\eps_n}(1+\eps_n\,v)|^{p_n}}{|x|^{b_np_n}}\;dx\\
&&=(\alpha_n+1)\int^{2\pi}_0\int^{+\infty}_0r^{2\frac{\alpha_n+\eps_n}{1-\eps_n}+1}\,\frac{\big(1+\eps_n\,v(r\cos\theta,r\sin\theta)\big)^{2/\eps_n}}{\big(1+r^{2(\alpha_n+1)}\big)^{\frac{2}{1-\eps_n}}}\;dr\,d\theta\\
&&=\int^{2\pi}_0\int^{+\infty}_0\frac{t^{\frac{1+\eps_n}{1-\eps_n}}}{(1+t^2)^{\frac{2}{1-\eps_n}}}\;\big(1+\eps_n\,v(\,t^{\frac{1}{1+\alpha_n}}\cos\theta,\;t^{\frac{1}{1+\alpha_n}}\sin\theta)\big)^{2/\eps_n}\;dt\,d\theta\,.
\end{eqnarray*}
We pass to the limit as $n \to + \infty$ and obtain:
\begin{eqnarray*}
\lim_{n\to +\infty}\frac1{\kappa_{\eps_n}}\int_{\R^2}\kern-3pt\frac{|u_{\eps_n}(1+\eps_n\,v)|^{p_n}}{|x|^{b_np_n}}\;dx&\kern -4pt=&\kern -4pt\frac1{\pi}\int^{2\pi}_0\kern-3pte^{2\,v(\cos\theta,\;\sin\theta)}\,d\theta\;\int^{+\infty}_0\kern-3pt\frac{t\;dt}{(1+t^2)^2}\\
&\kern -4pt=&\kern -4pt\frac1{2\pi}\int^{2\pi}_0e^{2\,v(\cos\theta,\;\sin\theta)}\,d\theta\;.
\end{eqnarray*}
Analogously,
\begin{eqnarray*}
&&\hspace*{-12pt}(\alpha_n+1)\int_{\R^2}\frac{u_{\eps_n}^{2/{\eps_n}}}{|x|^{2(a_n-\alpha_n)}}\,v\;dx\\
&&=(\alpha_n+1)\int_{\R^2}{|x|^{2\frac{\alpha_n+\eps_n}{1-\eps_n}}}\,\frac{v(x)}{(1+|x|^{2(1+\alpha_n)})^{\frac2{1-\eps_n}}}\;dx\\
&&=\int^{2 \pi}_0 \int^{+ \infty}_0 \frac{t^{\frac{1+\eps_n}{1 - \eps_n}}}{(1+t^2)^{{\frac{2}{1-\eps_n}}}}\;v(\,t^{\frac1{\alpha_n+1}}\cos\theta,\;t^{\frac{1}{\alpha_n+1}}\sin\theta)\;dt\,d\theta\,.
\end{eqnarray*}
By Lemma \ref{CP1}, we see that
\begin{eqnarray*}\frac1{\lambda_{\eps_n}}\int_{\R^2}\frac{|\,\nabla[u_{\eps_n}(1+{\eps_n}\,v)]\,|^2}{|x|^{2a_n}}\;dx&\kern -4pt=&\kern -4pt 1+\frac{\eps_n^2}{\lambda_{\eps_n}}\;\frac{8\,(\alpha_n+1)^2}{(1-\eps_n)^2}\int_{\R^2}\frac{u_{\eps_n}^{2/{\eps_n}}}{|x|^{2(a_n-\alpha_n)}}\,v\;dx\\
&&\qquad+\;O\(\frac{\eps_n}{1+\alpha_n}\)+O\(\frac{\eps_n\,a_n^2}{\lambda_{\eps_n}}\)\;,
\end{eqnarray*}
and so
\begin{eqnarray*}
\lim_{n\to +\infty}\(\frac1{\lambda_{\eps_n}}\int_{\R^2}\frac{|\,\nabla[u_{\eps_n}(1+{\eps_n}\,v)]\,|^2}{|x|^{2a_n}}\;dx\)^{\!1/\eps_n}\kern-10pt&\kern -4pt=&\kern -4pt e^{\frac{2}{\pi}\int^{2\pi}_0v(\cos\theta,\;\sin\theta)\,d\theta\int^{+\infty}_0\frac{t\,dt}{(1+t^2)^2}}\\
&\kern -4pt=&\kern -4pt e^{\frac 1\pi\;\int^{2\pi}_0v(\cos\theta,\;\sin\theta)\,d\theta}\,.
\end{eqnarray*}
Hence the validity of  \eqref{IneqAsympTruc} would imply that for all $v\in C_c^\infty(\R^2)$, there holds:
\[\label{ineqfalse}
\frac1{2\pi}\int^{2\pi}_0e^{2\,v(\cos\theta,\;\sin\theta)}\,d\theta\leq e^{\frac1{\pi}\,\int^{2\pi}_0v(\cos\theta,\;\sin\theta)\,d\theta}\;.
\]
But this is clearly impossible, since such an inequality is violated for instance by the function $v(x)=v(x_1,x_2)=x_1^2\,\eta(x)$, with $\eta$ a standard cut-off function such that $\eta(x)=1$ if $|x|\leq 1$, $\eta(x)=0$ if $|x|\geq 2$. \endproof

%%%%%%%%%%%%%%%%%%%%%%%%%%%%%%%%%%%%%%%%%%%%%%%%%%%%%%%%%%%%%%%%%%%%%%%%%%%%%%%
%%%%%%%%%%%%%%%%%%%%%%%%%%%%%%%%%%%%%%%%%%%%%%%%%%%%%%%%%%%%%%%%%%%%%%%%%%%%%%%
\section{Symmetry breaking}\label{sect3}
\newcommand{\w}{w}

This section is devoted to the proof of Theorem \ref{TTT1}, (i). We start by establishing Property (i'), which is weaker, but it follows as an easy consequence of the results of section~\ref{sect2}.

%%%%%%%%%%%%%%%%%%%%%%%%%%%%%%%%%%%%%%%%%%%%%%%%%%%%%%%%%%%%%%%%%%%%%%%%%%%%%%%
\subsection{Proof of Property (i')}

By Lemma \ref{Lemma1} and Kelvin transformation, we can reduce the proof to the case $a<0$. Let us argue by contradiction and assume that there exists $\eps_0\in (0,1)$, $a_n\to 0_-$ and $b_n$ such that $\eps_0<\frac{b_n}{a_n}<1$ and $u_{a_n, b_n}$ is radially symmetric. Set $\eps_n=b_n-a_n>0$ and define $\alpha_n$ such that $\alpha_n+1=-a_n\,(1-\eps_n)/\eps_n$. Notice that $\eps_n\to 0_+$ while $\alpha_n+1=a_n-a_n/(b_n-a_n)=a_n-(b_n/a_n-1)^{-1}>a_n+ (1-\eps_0)^{-1}$. Hence, $\liminf_{n\to +\infty}\alpha_n\ge\alpha_0=\eps_0/(1-\eps_0)>0$. But this is impossible since it contradicts Proposition \ref{CT2} in case $\liminf_{n\to +\infty}\alpha_n=+\infty$, or Propositions \ref{nottrue} and \ref{CT1} if $\limsup_{n\to +\infty}\alpha_n<+\infty$; and we conclude the proof of (i'). \endproof

%%%%%%%%%%%%%%%%%%%%%%%%%%%%%%%%%%%%%%%%%%%%%%%%%%%%%%%%%%%%%%%%%%%%%%%%%%%%%%%
\subsection{Proof of (i) of Theorem \ref{TTT1}}

It is well known (see \cite{Catrina-Wang-01}) that by means of the following Emden-Fowler transformations:
\bq\label{3.1}
t=\log|x|\;,\quad\theta=\frac{x}{|x|}\in S^{N-1}\,,\quad\w(t,\theta)=|x|^{\frac{N-2-2a}{2}}\,v(x)\;,
\eq
inequality~\eqref{HSN} for $u$ is equivalent to the Sobolev inequality for $w$ on $\R\times S^{N-1}$. Namely,
\[\label{3.3}
\|\w\|^2_{L^p(\R\times S^{N-1})}\leq\,C^N_{a,b}\;\Big[\|\nabla\w\|^2_{L^2(\R\times S^{N-1})}+{\textstyle \frac14(N-2-2a)^2}\|\w\|^2_{L^2(\R\times S^{N-1})}\Big]\,,
\]
for $\w\in H^1(\R\times S^{N-1})$, with $p=2\,N/[(N-2)+2\,(b-a)]$ and the same optimal constant $C^N_{a,b}$ as in~\eqref{HSN}. This inequality is consistent with the statement of Lemma~\ref{Lemma1}, as it makes sense for any $a\neq(N-2)/2$, independently of the sign of $N-2-2a$.

For $N=2$, the inequality holds for functions $\w=\w(t,\theta)$ defined over the two-dimensional cylinder $\mathcal C=\R\times S^1\approx(\R/2\pi\Z)$, {\sl i.e.,\/} such that $\w(t,\cdot)$ is $2\pi$-periodic for a.e. $t\in\R$. The inequality then takes the form
\bq\label{3.4}
\|\w\|^2_{L^p(\mathcal C)}\leq C_{a,a+2/p}\;\Big(\|\nabla\w\|^2_{L^2(\mathcal C)}+a^2\,\|\w\|^2_{L^2(\mathcal C)}\Big)\quad\forall\;\w\in H^1(\mathcal C)
\eq
for all $a\neq0$ and $p>2$. Here $C_{a,b}$ is the optimal constant in \eqref{HS} which enters in \eqref{3.4} with $b=a+2/p$.

For any $a\neq0$ and $p>2$, inequality~\eqref{3.4} is attained at an extremal function $\w_{a,p}\in H^1(\mathcal C)$ which satisfies
\bq\label{3.5}
\left\lbrace \begin{array}{ll}
-(\w_{tt}+\w_{\theta\theta})+a^2\,\w=\w^{p-1}\quad\mbox{in}\quad\R\times[-\pi,\pi]\;,\vspace*{12pt}\\
\w>0\;,\quad\w(t,\cdot)\quad\mbox{is $2\pi$-periodic}\quad\forall\;t\in\R\;,\end{array} \right.
\eq
and such that
\[\label{3.6}
\(C_{a,a+2/p}\)^{-1}=\|\w_{a,p}\|^{p-2}_{L^p(\mathcal C)}=\inf_{\w\in H^1(\mathcal C)\setminus\lbrace 0 \rbrace}\mathcal F(\w)\;,
\]
where the functional
\[\label{3.13}
\mathcal F(\w)=\frac{\|\nabla\w\|^2_{L^2(\mathcal C)}+a^2\,\|\w\|^2_{L^2(\mathcal C)}}{\|\w\|^2_{L^p(\mathcal C)}}
\]
is well defined on $H^1(\mathcal C)\setminus\lbrace 0\rbrace$. Moreover, according to \cite{Catrina-Wang-01}, we can further assume that
\bq\label{3.7-3.8}\left\{\begin{array}{l}
\w_{a,p}(t,\theta)=\w_{a,p}(-t,\theta)\quad\forall\;t\in\R\;,\quad\theta\in[-\pi,\pi)\;,\vspace{6pt}\\
\frac{\partial\w_{a,p}}{\partial t}\,(t,\theta)<0\quad\forall\;t>0\;,\quad\forall\;\theta\in[-\pi,\pi)\,,\vspace{6pt}\\
\max_{\R\times[-\pi,\pi)}\w_{a,p}=\w_{a,p}(0,0)\;.
\end{array}\right.\eq
This symmetry result is easy to establish for a minimizer, but the monotonicity requires more elaborate tools like the sliding method and we refer to~\cite{Catrina-Wang-01} for more details. For a solution of~\eqref{3.5} which does not depend on~$\theta$, the  conditions in  \eqref{3.7-3.8} allow to determine its value at $0$ simply by multiplying the ODE by $\w_t$ and integrating from $0$ to $\infty$. In fact, in this way, one deduces the relation: $a^2\,\w^2(0)/2=\w^p(0)/p$, which uniquely determines $\w(0)>0$.  In turn this yields to the following unique $\theta$-independent solution for~\eqref{3.5} and~\eqref{3.7-3.8}:
\[\label{3.9}
\w^*_{a,p}(t)=\big({\textstyle\frac{a^2\,p}{2}}\big)^{1/(p-2)}\,\left[\cosh\({\textstyle\frac{p-2}{2}\,a\,t}\)\right]^{-2/(p-2)}\,,
\]
as a consequence of the classification result in~\cite{Catrina-Wang-01}. Such a solution is an extremal for \eqref{3.4} on the set of functions which are independent of the $\theta$-variable, and
\[\label{3.10}
\|\w^*_{a,p}\|^{\,p-2}_{L^p(\R)}=\inf_{f\in H^1(\R)\setminus\lbrace 0\rbrace}\mathcal F^*(f)\quad\mbox{with}\quad\mathcal F^*(f)=\frac{\|f'\|^2_{L^2(\R)}+a^2\,\|f\|^2_{L^2(\R)}}{\|f\|^2_{L^p(\R)}}\;.
\]
For simplicity, we will also write $\mathcal F(f)=(\pi)^{1-2/p}\,\mathcal F^*(f)$ for all functions $f$ which are independent of $\theta$. As a useful consequence of the above considerations, we have the following result.
%------------------------------------------------------------------------------
\begin{lemma}\label{Remark9} Let $p>2$. For any $a\neq0$,
\[
\(C_{a,a+2/p}\)^{-\frac p{p-2}}=\|\w_{a,p}\|^p_{L^p(\mathcal C)}\leq\|\w^*_{a,p}\|^p_{L^p(\mathcal C)}= 4\,\pi\,(2\,a)^\frac p{p-2}\,(a\,p)^\frac 2{p-2}\,c_p
\]
where $c_p$ is an increasing function of $p$ such that
\bq\label{Behavior-Cp}\begin{array}{lll}
c_p\to 0\quad&\mbox{as}\quad&p\to2_+\;,\vspace{6pt}\\
c_p\to \frac 12\quad&\mbox{as}\quad&p\to+\infty\;.
\end{array}\eq
As a consequence, if $a=a(p)$ is such that $\lim_{p\to\infty}a(p)\,p=2\,(\alpha+1)$, then
\bq\label{Bidule-a}
\lim_{p\to\infty} \,p\int_\mathcal C|w_{a(p),p}^*|^p\,dx=8\,(\alpha+1)\;.
\eq
\end{lemma}
%------------------------------------------------------------------------------
\proof Observe that
\[
\|\w_{a,p}\|^p_{L^p(\mathcal C)}\!=\(C_{a,a+2/p}\)^{-\frac p{p-2}}\!=\(\mathcal F(\w_{a,p})\)^\frac p{p-2}\leq\(\mathcal F(\w^*_{a,p})\)^\frac p{p-2}\!=\|\w^*_{a,p}\|^p_{L^p(\mathcal C)}\,.
\]
On the other hand,
\begin{eqnarray*}\label{GHG}
\|\w^*_{a,p}\|^p_{L^p(\mathcal C)}&=&2\pi\,\big({\textstyle\frac{a^2\,p}{2}}\big)^{\frac p{p-2}}\int_{-\infty}^\infty\left[\cosh\Big({\textstyle\frac{a\,(p-2)}2\,t}\Big)\right]^{-\frac{2p}{p-2}}\,dt\\
&=&4\pi\,\big({\textstyle\frac{a^2\,p}{2}}\big)^{\frac p{p-2}}\int_0^\infty\frac{2^{\frac{2p}{p-2}}\,e^{-a\,p\,t}}{\big(1+e^{-a\,(p-2)\,t}\big)^{\frac{2p}{p-2}}}\;dt\\
&=&4\pi\,\big({\textstyle\frac{a^2\,p}{2}}\big)^{\frac p{p-2}}\;{\textstyle\frac{2^{\frac{2p}{p-2}}}{a\,p}}\int_0^1\frac{ds}{\big(1+s^{(p-2)/p}\big)^{\frac{2p}{p-2}}}\;.
\end{eqnarray*}
Hence by setting:
$$
c_p=\int_0^1\frac{ds}{\big(1+s^{(p-2)/p}\big)^{\frac{2p}{p-2}}}\;,
$$
we easily check \eqref{Behavior-Cp} and the fact that $c_p$ is monotonically increasing in $p$. The limiting behavior of $c_p$ stated in \eqref{Behavior-Cp} is a direct consequence of Lebesgue's dominated convergence theorem. \endproof
We can now reformulate Theorem \ref{TTT1} in the cylinder $\mathcal C$, in terms of $\w$, as follows.
%------------------------------------------------------------------------------
\begin{theorem}\label{Theorem2} Let $a\neq0$ and $p>2$.
\begin{itemize}
\item[(i)] If $|a|\,p>2\sqrt{1+a^2}$, then $\mathcal F(\w_{a,p})<\mathcal F(\w^*_{a,p})$.
\item[(ii)] For every $\eps>0$, there exists $\delta>0$ such that, if $0<|a|<\delta$ and $|a|\,p<2-\eps$, then $\mathcal F(\w_{a,p})=\mathcal F(\w^*_{a,p})$.
\end{itemize}
\end{theorem}
%------------------------------------------------------------------------------
Part (ii) of Theorem \ref{Theorem2} will be proved in the next section. Concerning part (i), we define the quadratic form
\[\label{3.14}
Q(\psi)=\|\nabla\psi\|^2_{L^2(\mathcal C)}+a^2\,\|\psi\|^2_{L^2(\mathcal C)}-(p-1)\int_\mathcal C|\w^*_{a,p}|^{p-2}\,|\psi|^2\,dx
\]
on $H^1(\mathcal C)$. In fact, property (i) is a consequence of the following result, inspired by~\cite{Catrina-Wang-01,Felli-Schneider-03} (at least for the case $a<0$):%------------------------------------------------------------------------------
\begin{proposition}\label{Proposition3.1} Let $a\neq0$ and $p>2$. Then
\[
\inf_{\begin{array}{c}\psi\in H^1(\mathcal C)\\ \int^\pi_{-\pi}\psi(t,\theta)\,d\theta=0,\;t\in\R\;\mbox{a.e.}\end{array}}\hspace*{-24pt}\frac{Q(\psi)}{\|\psi\|_{L^2(\mathcal C)}^2}\;=\;a^2+1-\(\frac{a\,p}2\)^2
\]
is achieved by
\[
\psi(t,\theta)=\big(\cosh((\alpha+1)\,t)\big)^{-p/(p-2)}\,,\quad\mbox{with}\quad\alpha=(p-2)\,a/2-1\;.
\]
In particular, if $|a|\,p>2\sqrt{1+ a^2}$\,, then $\w^*_{a,p}$ is a critical point for $\mathcal F$ of saddle-type. \end{proposition}
%------------------------------------------------------------------------------
\proof Since $\w^*_{a,p}$ is a local minimum for $\mathcal F$ when restricted to the set of functions independent of $\theta$, to search for negative directions of the Hessian of $\mathcal F$ around $\w^*_{a,p}$, we have to analyze the quadratic form $Q(\psi)$ on the space of functions $\psi\in H^1(\mathcal C)$ such that $\int^\pi_{-\pi}\psi(t,\theta)\,d\theta=0$ for a.e. $t\in\R$. To this purpose, we use the Fourier expansion of $\psi$,
\begin{eqnarray*}\label{3.15}
&&\hspace*{-20pt}\psi(t,\theta)=\sum_{k\neq0}f_k(t)\;\frac{e^{ik\theta}}{\sqrt{2\pi}}\;,\quad f_{-k}(t)=\overline{f_k}(t)\;,\\
&&\hspace*{-20pt}Q(\psi)=2\sum^{+\infty}_{k=1}\Big(\|f'_k\|^2_{L^2(\R)}\!+\!(a^2\!+\!k^2)\,\|f_k\|^2_{L^2(\R)}\!-\!(p\!-\!1)\int_\R|\w^*_{a,p}|^{p-2}\,|f_k|^2\,dt\Big)\,.
\end{eqnarray*}
Hence we obtain a negative direction for $Q$ if and only if
\[\label{3.16}
\mu^1_{a,p}=\!\!\!\inf_{f\in H^1(\R)\setminus\lbrace 0\rbrace}\!\!\!\!\frac{\|f'\|^2_{L^2(\R)}\!\!\!+\!(a^2+1)\|f\|^2_{L^2(\R)}\!\!\!-(p-1)\int_{\R}|\w^*_{a,p}|^{p-2}\,|f|^2\,dt}{\|f\|^2_{L^2(\R)}}<0\;.
\]
Setting $1+\alpha=(p-2)\,a/2$ and $\beta=a^2\,p\,(p-1)/2=2\,(1+\alpha)^2\,p\,(p-1)/(p-2)^2>0$, the question is reduced to the eigenvalue problem
\[\label{3.18}
-f''-\frac{\beta\,f}{\big(\cosh((\alpha+1)\,t)\big)^2}=\lambda\,f\;.
\]
in $H^1(\R)$. The eigenfunction $f_1(t)=\big(\cosh((\alpha+1)\,t)\big)^{-p/(p-2)}$ corresponds to the first eigenvalue $\label{38ter}\lambda_1=-(a\,p/2)^2$. See \cite{Landau-Lifschitz-67, Felli-Schneider-03} for a discussion of the above eigenvalue problem. Hence $\mu^1_{a,p}=1+a^2-(a\,p/2)^2$, and the proof is completed. \endproof

%%%%%%%%%%%%%%%%%%%%%%%%%%%%%%%%%%%%%%%%%%%%%%%%%%%%%%%%%%%%%%%%%%%%%%%%%%%%%%%
%%%%%%%%%%%%%%%%%%%%%%%%%%%%%%%%%%%%%%%%%%%%%%%%%%%%%%%%%%%%%%%%%%%%%%%%%%%%%%%
\section{A symmetry result}\label{sect4}

The section is devoted to the proof of part (ii) of Theorem~\ref{Theorem2}.

Without loss of generality, by Lemma~\ref{Lemma1}, we can restrict our analysis to the case $a>0$.

%%%%%%%%%%%%%%%%%%%%%%%%%%%%%%%%%%%%%%%%%%%%%%%%%%%%%%%%%%%%%%%%%%%%%%%%%%%%%%%
\subsection{Pohozaev's identity}
%------------------------------------------------------------------------------
\begin{lemma}\label{Lemma1bis} If $\w\in H^1(\mathcal C)$ satisfies \eqref{3.5}, then for all $t\in\R$, $w=w(t,\theta)$ satisfies the identity
\[\label{e11}
\int^\pi_{-\pi}\Big(\frac{\partial\w}{\partial\theta}\Big)^2\;d\theta = \int^\pi_{-\pi}\Big(\frac{\partial\w}{\partial t}\Big)^2\;d\theta - a^2\int^\pi_{-\pi}\w^2\;d\theta+\frac{2}{p}\int^\pi_{-\pi}\w^p\;d\theta\,.
\]\end{lemma}
%------------------------------------------------------------------------------
\proof Multiply the equation in \eqref{3.5} by $\frac{\partial\w}{\partial t}$ and integrate over $[- \pi, \pi]$ to obtain:
$$
\int^\pi_{-\pi}\left(-\frac{\partial^2\w}{\partial t^2}\,\frac{\partial\w}{\partial t}-\frac{\partial^2\w}{\partial\theta^2}\,\frac{\partial\w}{\partial t}+a^2\,\frac{\partial\w}{\partial t}\,\w\right)\,d\theta=\int^\pi_{-\pi}\w^{p-1}\,\frac{\partial\w}{\partial t}\,d\theta\;,
$$
that is
$$
\int^\pi_{-\pi}\left\{{\textstyle -\frac{\partial}{\partial\theta}\(\frac{\partial\w}{\partial\theta}\frac{\partial\w}{\partial t}\)+\frac 12\,\frac{d}{dt}\left[\(\frac{\partial\w}{\partial\theta}\)^2-\(\frac{\partial\w}{\partial t}\)^2+a^2\,\w^2\right]}\right\}d\theta=\frac1p\int^\pi_{-\pi}\frac{d\,(\w^p)}{dt}\,d\theta\;.
$$
Since $\int^\pi_{-\pi}\frac{\partial}{\partial\theta}\(\frac{\partial\w}{\partial\theta}\frac{\partial\w}{\partial t}\)\,d\theta=0$, we get
$$
\frac{d}{dt}\int^\pi_{-\pi}\left[\(\frac{\partial\w}{\partial t}\)^2-\(\frac{\partial\w}{\partial\theta}\)^2-a^2\,\w^2+\frac{2}{p}\,\w^p\right]d\theta= 0
$$
for all $t \in \R$. Hence as a function of $t$, the above integral must be a constant. Since it is also integrable over $\R$, then it must vanish identically. \endproof

%%%%%%%%%%%%%%%%%%%%%%%%%%%%%%%%%%%%%%%%%%%%%%%%%%%%%%%%%%%%%%%%%%%%%%%%%%%%%%%
\subsection{Proof of Theorem \ref{Theorem2}, (ii)}

We argue by contradiction and suppose that there exists $\eps_0 \in (0,1)$ and, for all $n\in\N$, $a_n>0$, $p_n>2$, such that:
\bq\label{e12}
\lim_{n\to+\infty}a_n = 0\;,\quad a_n\,p_n<2-\eps_0\quad\mbox{and}\quad\mathcal F(\w_{a_n,\,p_n})<\mathcal F(\w^*_{a_n,\,p_n})\;.
\eq
For simplicity, set
\[\label{e13}
\w_n=\w_{a_n,\,p_n}\quad\mbox{and}\quad\w^*_n=\w^*_{a_n,\,p_n}\,,
\]
and recall that we can assume
$$
\w_n(t,\theta)=\w_n(-t,\theta)\;,\quad\frac{\partial\w_n}{\partial t}(t,\theta)<0\quad\forall\;t>0\quad\mbox{and}\quad\w_n(0,0)=\max_\mathcal C\,\w_n\;.
$$
Notice in particular that $\frac{\partial\w_n}{\partial t}(0,\theta)=0$ for any $\theta\in[-\pi,\pi]$. If we apply Lemma~\ref{Lemma1bis} to $\w=\w_n$ and $t=0$, we obtain
$$
\frac{p^2_n\,a^2_n}2\int^\pi_{-\pi}\w^2_n(0,\theta)\,d\theta\leq p_n\!\int^\pi_{-\pi}\w_n^{p_n}(0,\theta)\,d\theta\leq p_n\,\|\w_n\|^{p_n-2}_{L^\infty(\mathcal C)}\!\int^\pi_{-\pi}\w^2_n(0,\theta)\,d\theta\;,
$$
and deduce that
$$
p_n\,\|\w_n\|^{p_n-2}_{L^\infty(\mathcal C)}\geq\frac 12\,p^2_n\,a^2_n\;.
$$
%------------------------------------------------------------------------------
\begin{lemma}\label{TL1}
\[\label{e15}\liminf_{n\to+\infty}p_n\,\|\w_n\|^{p_n-2}_{L^\infty(\mathcal C)}\geq 1\;.
\]
\end{lemma}
%------------------------------------------------------------------------------
\proof We can write $\w_n(t,\theta)=\varphi_n(t)+\psi_n(t,\theta)$ with
\[\varphi_n(t)=\frac1{2\pi}\int^\pi_{-\pi}\w_n(t,\theta)\,d\theta\;,\quad\kern -3pt\int^\pi_{-\pi}\psi_n(t,\theta)\,d\theta=0\quad\kern -3pt\mbox{a.e.}\;t\in\R\quad\kern -3pt\mbox{and}\quad\kern -3pt\psi_n\neq 0\;.
\]
Multiplying \eqref{3.5} by $\psi_n$ and using the fact that $\int^\pi_{-\pi}\psi_n(t,\theta)\,d\theta=0$ for any $t\in\R$, we find
\begin{eqnarray*}
&&\left\|\frac{\partial\psi_n}{\partial t}\right\|^2_{L^2(\mathcal C)}+\;\left\|\frac{\partial\psi_n}{\partial\theta}\right\|^2_{L^2(\mathcal C)}+\;a^2_n\,\|\psi_n\|^2_{L^2}\\
&&=\int_\mathcal C\w^{p_n-1}_n\,\psi_n\;dt\,d\theta\\
&&=\int_\mathcal C\w^{p_n-1}_n\,\psi_n\;dt\,d\theta-\int_\mathcal C\varphi_n^{p_n-1}\,\psi_n\;dt\,d\theta\\
&&=(p_n-1)\int^1_0\left\{\int_\mathcal C\big|s\,\varphi_n+(1-s)\,\w_n\big|^{p_n-2}\,|\psi_n|^2\;dt\,d\theta\right\}\,ds\\
&&\label{e17} \leq(p_n-1)\,\|\w_n\|^{p_n-2}_{L^\infty(\mathcal C)}\int_\mathcal C|\psi_n|^2\;dt\,d\theta\;.
\end{eqnarray*}
By Poincar\'e's inequality, we know that $\|\psi_n\|^2_{L^2(\mathcal C)}\leq\big\|\frac{\partial\psi_n}{\partial\theta}\big\|^2_{L^2(\mathcal C)}$, and this proves the claim. \endproof

Next we introduce the new parameters:
\[\label{e19}
\eps_n=\frac{2}{p_n}\quad\mbox{and}\quad\alpha_n=-1+(1-\eps_n)\,\frac{a_n}{\eps_n}=-1+\frac 12\,(1-\eps_n)\,a_n\,p_n\;.
\]
%------------------------------------------------------------------------------
\begin{lemma}\label{TL4} Up to a subsequence we have:\[\label{e20}
\lim_{n\to+\infty}\alpha_n=\alpha\in[-1,0)\;,
\]
and $\lim_{n\to+\infty} p_n=+\infty$, or equivalently,
\[
\lim_{n\to+\infty}\eps_n=0\;.
\]
\end{lemma}
%------------------------------------------------------------------------------
\proof From the condition: $a_n\,p_n<2-\eps_0$, we deduce that $\alpha_n+1\leq(1-\eps_n)\,(1-\eps_0/2)$. Thus, along a subsequence, we can assume that $\alpha_n$ converges to some $\alpha\in[-1,0)$ and $\lim_{n\to+\infty}p_n\in[2,\infty]$.

To rule out the possibility that $\lim_{n\to+\infty}p_n=\bar p\in[2,\infty)$, notice that if this would be the case, then by Lemma~\ref{Remark9},
\[
\lim_{n\to+\infty}\|\w_n\|_{L^{p_n}(\mathcal C)}=0\;.
\]
By applying local elliptic estimates in a neighborhood of the origin $(0,0)$ then we would deduce that $\lim_{n\to+\infty}\|\w_n\|_{L^{\infty}(\mathcal C)}=\lim_{n\to+\infty}\w_n(0,0)=0$, in contradiction with Lemma~\ref{TL1}. \endproof

%------------------------------------------------------------------------------
\begin{corollary}\label{CTL1}
\[
\liminf_{n\to+\infty}w_n(0,0)\geq 1\;.
\]
\end{corollary}
%------------------------------------------------------------------------------
\proof If by contradiction we assume that $\liminf_{n\to+\infty}w_n(0,0)<1$, then $\liminf_{n\to+\infty}p_n\,\|\w_n\|_{L^\infty(\mathcal C)}^{p_n-2}=0$, and again this is impossible by Lemma~\ref{TL1}. \endproof

%------------------------------------------------------------------------------
\begin{lemma}\label{TL2}
\[\label{e23}
\limsup_{n\to+\infty}\,p_n\,\|\w_n\|^{p_n-2}_{L^\infty(\mathcal C)}<+\infty\;.
\]
\end{lemma}
%------------------------------------------------------------------------------
\proof Argue by contradiction, and assume that, along a subsequence, $\delta_n=\big(p_n\,\|\w_n\|^{p_n-2}_{L^\infty(\mathcal C)}\big)^{-1/2}$ converges to $0$ as $n\to+\infty$. We consider the function
\[\label{e24}
\W_n(t,\theta)=p_n\,\(\frac{\w_n(\delta_n\,t,\,\delta_n\,\theta)}{\w_n(0,0)}-1\)
\]
defined in $\mathcal C_n=\R\times[-\pi/\delta_n,\,\pi/\delta_n]$, which satisfies
\bq\label{e25}\left\lbrace\begin{array}{ll}
-\Delta\W_n=\(1+\frac{\W_n}{p_n}\)^{p_n-1}\kern -8pt-\,a^2_n\,p_n\,\delta^2_n\,\(1+\frac{\W_n}{p_n}\)\quad\mbox{in}\quad\mathcal C_n\;,\vspace*{12pt}\\
\W_n\leq 0=\W_n(0,0)\;.\end{array} \right.
\eq
Furthermore, by Lemma~\ref{Remark9}, we find\
\[
\label{e26}
\int_{\mathcal C_n}\(1+\frac{\W_n}{p_n}\)^{p_n}\,dx=\frac{p_n}{\w_n(0,0)^2}\,\int_\mathcal C\w_n^{p_n}\,dx\leq \frac 1{\w_n(0,0)^2}\;p_n\int_\mathcal C|\w_n^*|^{p_n}\,dx\;.
\]
Recalling that  $\liminf_{n\to+\infty}\w_n(0,0)\geq 1$ and $\lim_{n\to+\infty}a_n\,p_n=2\,(1+\alpha)$ by~\eqref{Bidule-a}, we can pass to the limit above and by virtue of \eqref{Behavior-Cp}-\eqref{Bidule-a}, conclude:
\[\label{new-est}
\lim_{n\to+\infty}\|1+\W_n/p_n\|_{L^{p_n}(\mathcal C_n)}^{p_n}\leq 8\pi\,(1+\alpha)\;.
\]

Since the right hand side in \eqref{e25} is uniformly bounded in $L^\infty_{\rm loc} (\R^2)$, we can use Harnack's inequality (see for instance \cite{Brezis-Merle-91, Tarantello-07} in similar cases) to deduce that $\W_n$ is uniformly bounded in $L^\infty_{\rm loc}$. Hence, by elliptic regularity theory, $\W_n$ is uniformly bounded in $C^{2,\alpha}_{\rm loc}$. So we can find a subsequence along which $\W_n$ converges pointwise (uniformly on every compact set in $\R^2$) to a function $\W$ which satisfies
\begin{equation}\label{e27}
-\Delta\W=e^{\W}\quad\mbox{in}\quad\R^2\,.
\end{equation}
Furthermore, by Fatou's Lemma,
\[\label{e27'}
\int_{\R^2}e^{\W}\,dx\leq\lim_{n\to+\infty}\;\int_{\mathcal C_n}\(1+\frac{\W_n}{p_n}\)^{p_n}dx\leq 8\pi\,(1+\alpha)<8\pi\;,
\]
as $\alpha\in[-1,0)$. But this is impossible, since according to \cite{MR1121147}, every solution~$\W$ of \eqref{e27} with $e^{\W}\in L^1(\R^2)$, must satisfy $\int_{\R^2}e^{\W}\,dx=8\pi$ (also see \cite{Chou-Wan-94,MR1372247}). \endproof

%------------------------------------------------------------------------------
\begin{corollary}\label{TL3} For a subsequence of $\|\w_n\|_{L^\infty(\mathcal C)}=w_n(0,0)$ (denoted the same way) we have:
\begin{eqnarray*}
&&\lim_{n\to+\infty}w_n(0,0)=1\;,\\
&&\lim_{n\to+\infty}\big[w_n(0,0)\big]^{p_n}=0\;,\\
&&\lim_{n\to+\infty}p_n\,\big[w_n(0,0)\big]^{p_n-2}=\mu\in [1,+\infty)\;.
\end{eqnarray*}
\end{corollary}
%------------------------------------------------------------------------------
\proof The existence of a limit $\mu\geq 1$ is just a consequence of Lemmata~\ref{TL1} and \ref{TL2}. Furthermore by Lemma~\ref{TL4}, $p_n=2/\eps_n\to+\infty$ as $n\to+\infty$, which proves that $[w_n(0,0)]^{p_n}$ converges to $0$. Finally, according to Corollary~\ref{CTL1}, $\liminf_{n\to+\infty}w_n(0,0)\geq 1$ and if this limit were not $1$, we would get a contradiction to the existence of $\mu$. \endproof

Define the function
\[
V_n(t,\theta)=p_n\(\frac{\w_n(t,\theta)}{\w_n(0,0)}-1\)\quad\forall\;(t,\theta)\in\mathcal C\;.
\]
It satisfies:
\begin{eqnarray*}
&&-\Delta V_n=p_n\,\big(\w_n(0,0)\big)^{p_n-2}\(1+\frac{V_n}{p_n}\)^{p_n-1}\kern -8pt-a^2_n\,p_n\(1+\frac{V_n}{p_n}\)\quad\mbox{in}\quad\mathcal C\,,\\
&&\label{e30}V_n\leq 0=V_n(0,0)\;,\qquad V_n(t,\cdot)\quad \mbox{is $2\pi$-periodic}\;.
\end{eqnarray*}
We also observe that
\[\label{Bidule-b}
p_n\,\big(\w_n(0,0)\big)^{p_n}\int_\mathcal C\Big(1+\frac{V_n}{p_n}\Big)^{p_n}\,dx=p_n\,\int_\mathcal C|w_n|^{p_n}\,dx\leq p_n\,\int_\mathcal C|w_n^*|^{p_n}\,dx
\]
and by \eqref{Bidule-a}, $\lim_{n\to\infty}p_n\int_\mathcal C|w_n^*|^{p_n}\,dx=8\,\pi\,(\alpha+1)$. In particular, by Corollary~\ref{TL3}, we obtain
\[
\lim_{n\to+\infty}\;p_n\,\big(\w_n(0,0)\big)^{p_n-2}\int_\mathcal C\Big(1+\frac{V_n}{p_n}\Big)^{p_n}\,dx\,\leq\,8\pi\,(1+\alpha)\;.
\]
%------------------------------------------------------------------------------
\begin{lemma}\label{TL5} Up to a subsequence, $V_n$ converges to a function $V$ pointwise and $C^2$-uniformly on any compact set in $\R\times[-\pi,\pi]$. Furthermore $V$ satisfies:
\bq\label{e32}\left\lbrace\begin{array}{lll}
-\Delta V=\mu\,e^{V}\quad\mbox{in}\quad\mathcal C\;,\\ \\
\max_{\mathcal C}V\leq 0=V(0,0)\;,\quad V(t,\cdot)\quad\mbox{is $2\pi$-periodic}\quad\forall\;t\in\R\;,\\ \\
\mu\int_\mathcal C\,e^{V}\,dx\,\leq\,8\pi\,(1+\alpha)\;,\\
\end{array}\right.
\eq
\[\label{e33}
V(t,\theta)=V(-t,\theta)\;,\quad\frac{\partial V}{\partial t}(t,\theta)<0\quad\forall\;t>0\;,\quad\forall\;\theta\in[-\pi,\pi]\;,
\]
and
\bq\label{e35}
\int^\pi_{-\pi}\Big(\frac{\partial V}{\partial\theta}\Big)^2\,d\theta=\int^\pi_{-\pi}\Big(\frac{\partial V}{\partial t}\Big)^2\,d\theta-8\pi\,(1+\alpha)^2+2\mu\int^\pi_{-\pi}e^{V}\,d\theta\quad\forall\;t\in\R\;.
\eq
\end{lemma}
%------------------------------------------------------------------------------
\proof Since $-\Delta V_n$ is uniformly bounded in $L^\infty_{\rm loc} (\R^2)$, by Harnack's inequality, we see that $V_n$ is uniformly bounded in $L^\infty_{\rm loc}$. Hence, by elliptic regularity theory, $V_n$ is uniformly bounded in $C^{2,\alpha}_{\rm loc}$. Therefore, up to a subsequence, $V_n$ converges pointwise, and uniformly on every compact set in $\mathcal C$, to a function~$V$ which satisfies \eqref{e32} with $0\leq 1+\alpha< 1$, and also inherits the symmetric properties of $V_n$. To obtain~\eqref{e35} observe first that the result of Lemma~\ref{Lemma1bis} can be rewritten as follows,
\[
\int^\pi_{-\pi}\kern -3pt\Big(\frac{\partial V_n}{\partial\theta}\Big)^2d\theta=\int^\pi_{-\pi}\kern -3pt\Big(\frac{\partial V_n}{\partial t}\Big)^2d\theta - \frac{a^2_n\,p^2_n}{\w^2_n(0,0)}\int^{\pi}_{-\pi}\kern -8pt|\w_n|^2\,d\theta+\frac{2\,p_n}{\w^2_n(0,0)}\int_{-\pi}^\pi\kern -8pt|\w_n|^{p_n}\,d\theta\,,
\]
for any $t\in\R$, and that $w_n$ converges uniformly to $1$ on any compact set in $\R\times[-\pi,\pi]$. Hence by means of Lemma~\ref{TL4} and Corollary~\ref{TL3}, we can pass to the limit in the above identity and deduce~\eqref{e35}. \endproof
%------------------------------------------------------------------------------
\begin{lemma}\label{Lem:Chou-Wan} The following estimates hold:
\begin{eqnarray*}
&&\lim_{n\to+\infty}p_n\,\Big(\|\w_n\|^{p_n}_{L^{p_n}(\mathcal C)}-\|\w^*_n\|^{p_n}_{L^{p_n}(\mathcal C)}\Big)=0\;,\\
&&\int_\mathcal Ce^{V}\,dx=\lim_{n\to+\infty}\;\int_{\mathcal C_n}\(1+\frac{V_n}{p_n}\)^{p_n}\,dx=\frac{4\pi}{\alpha+1}\;.
\end{eqnarray*}
Moreover,
\[
\mu=2\,(\alpha+1)^2\,,
\]
and $V$ takes the form
\bq\label{ExtremalCylinder}
V(t)=-2\log\big[\cosh((\alpha+1)\,t)\big]\;.
\eq
\end{lemma}
%------------------------------------------------------------------------------
\proof In order to identify the given solution of \eqref{e32}, we consider the function $\varphi$ expressed in polar coordinates as follows:
\[
\varphi(r,\theta)=V\(-\log r\,,\,\theta\)-2\log r+\log\mu\quad\forall\;r>0\;,\quad\forall\;\theta\in[-\pi,\pi]\;.
\]
By straightforward calculations we see that $\varphi$ satisfies:
\begin{eqnarray*}
&&-\Delta\varphi=-\frac{1}{r^2}\,\(V_{tt}+V_{\theta\theta}\)\,\(-\log r\,,\,\theta\)=e^\varphi\quad\mbox{in}\;\R^2\backslash\lbrace0\rbrace\;,\\
&&\int_{\R^2}e^\varphi\,dx\,\leq\,8\pi\,(1+\alpha)\;,\end{eqnarray*}
and
\bq\label{e36}\varphi\(r^{-1}\,,\,\theta\)=\varphi(r,\theta)+4\,\log r\quad\forall\;r>0\;,\quad\forall\;\theta\in[-\pi,\pi]\;.\eq
A classification result of Chou and Wan (see \cite[Theorem 3, 1.]{Chou-Wan-94} and \cite{MR1372247}) concerning  solutions of Liouville equations on the punctured disk allows us to conclude that (in complex notations):
\[
\varphi(z)=\log\left[\frac{8\,|f'(z)|^2}{\big(1+|f(z)|^2\big)^2}\right]\,,
\]
with $f$ locally univalent in $\C\,\backslash \lbrace 0 \rbrace$, possibly multivalued and,
\begin{enumerate}
\item[(i)] either $f(z)=z^\gamma\,g(z)$,
\item[(ii)] or $f(z)=\phi(\sqrt z)$ and $\phi(z)\,\phi(-z)=1$,
\end{enumerate}
where $g$ and $\phi$ are holomorphic in $\C\,\backslash \lbrace 0 \rbrace$. Since the case (ii) implies that~$\phi$ must admit an essential singularity either at the origin or at infinity, this can be excluded in account of the integrability condition of $e^\varphi$.

On the other hand, in case (i), if we take into account the fact that $f'\neq 0$ for any $z\neq 0$, and the integrability of $e^\varphi$, we can allow only the choice:
$$
f(z)=a\,\big(z^{\beta+1}-b)\;,
$$
with $\beta\in\R$, $a$, $b\in\C$ and $b\neq 0$ only if $\beta+1\in\N$ (as otherwise $\varphi$ would be multivalued). For the corresponding solution $\varphi$ we find:
$$
\varphi(z)=\log\left[\frac{8\,\lambda\,(\beta+1)^2\,|z|^{2\beta}}{(1+\lambda\,|z^{\beta+1}-b|^2)^2}\right]\,,\quad\mbox{with}\quad\lambda=|a|^2\,.
$$
The symmetry property \eqref{e36} implies that
$$
\varphi\(\frac{z}{|z|^2}\)=\varphi(z)+4\log|z|\;,
$$
and so, necessarily $b=0$ and $\lambda=1$. Hence,
$$
\varphi(z)=\varphi(r)=\log\left[\frac{8\,(\beta+1)^2\,r^{2\beta}}{\big(1+r^{2(\beta+1)}\big)^2}\right]\,.
$$
By direct calculation, we get
$$\int_{\R^2}e^\varphi\,dx=8\pi\,(1+\beta)\,\leq\,8\pi\,(1+\alpha)\;.$$
In other words, $-1<\beta\leq\alpha<0$. As a consequence, we find that $V=V(t)$ is given by
$$
V(t)=\varphi(e^{-t})-2t-\log\mu=\log\left[\frac{2\,(\beta+1)^2}{\mu\,\big(\cosh((\beta+1)\,t)\big)^2}\right]\;,
$$
with $-1<\beta\leq\alpha<0$. The condition $V(0)=0$ implies $\mu=2\,(\beta+1)^2$.

\noindent On the other hand, from \eqref{e35} we also have:
$$
\Big(\frac{\partial V}{\partial t}\Big)^2=4\,(1+\alpha)^2-\frac{4\,(\beta+1)^2}{\big(\cosh((\beta+1)\,t)\big)^2}\;,
$$
that gives:
$$
4\,(\beta+1)^2\frac{\big(\sinh((\beta+1)\,t)\big)^2}{\big(\cosh((\beta+1)\,t)\big)^2}=4\,(1+\alpha)^2-\frac{4\,(\beta+1)^2}{\big(\cosh((\beta+1)\,t)\big)^2}\;,
$$
and we get $\beta=\alpha$. Therefore \eqref{ExtremalCylinder} is established and necessarily
\[
\lim_{n\to\infty}p_n\,\big(\w_n(0,0)\big)^{p_n}\int_\mathcal C\Big(1+\frac{V_n}{p_n}\Big)^{p_n}\,dx=2(\alpha+1)^2\int_{\mathcal C}e^V\,dx=8\,\pi\,(\alpha+1)\;.
\]
Thus, by recalling \eqref{Bidule-a}, we complete the proof. \endproof

Define
\[
r_n=\sup_\mathcal C\,\bigg|\(\frac{\w_n}{\w_n(0,0)}\)^{p_n-2}-\,e^{V}\bigg|\;.
\]
%------------------------------------------------------------------------------
\begin{lemma}\label{Proposition1} With the above notations, $\lim_{n\to+\infty} r_n=0$. \end{lemma}
%------------------------------------------------------------------------------
\proof Fix $\eps > 0$ and choose $R_\eps>0$ sufficiently large so that
$$
e^{V(R_\eps)}=\frac{1}{\big(\cosh((\alpha+1)\,R_\eps)\big)^2}\,<\;\frac{\eps}{4}\,.
$$
Furthermore, $\(\w_n(t,\theta)/\w_n(0,0)\)^{p_n-2}=(1+V_n/p_n)^{p_n-2}$ converges to $e^{V}$ uniformly on any compact set in $\R\times[-\pi,\pi]$, and so we can find $n_\eps\in\N$ such that for all $n\geq n_\eps$,
$$
\sup_{\scriptstyle |t|\leq R_\eps\,,\;|\theta|\leq\pi}\Big|\Big(\frac{\w_n(t,\theta)}{\w_n(0,0)}\Big)^{p_n-2}-e^{V}\Big|\,<\;\frac{\eps}{4}\,.$$
Thus, recalling that $\(\w_n(t,\theta)/\w_n(0,0)\)^{p_n-2}$ and $e^V$ are even in $t$ and monotone decreasing in $t>0$ by Lemma~\ref{TL5}, for $n\geq n_\eps$ we find the estimate
\[
r_n\leq\underbrace{\sup_{\scriptstyle|t|\leq R_\eps\,,\;|\theta|\leq\pi}\Big|\Big(\frac{\w_n(t,\theta)}{\w_n(0,0)}\Big)^{p_n-2}\!-e^{V}\Big|}_{<\eps/4}+\underbrace{\sup_{|t|\geq R_\eps}\Big(\frac{\w_n(t,\theta)}{\w_n(0,0)}\Big)^{p_n-2}}_{e^{V(R_\eps)}+\eps/4<\eps/2}+\underbrace{\sup_{|t|\geq R_\eps}e^{V}}_{\eps/4}\;,
\]
which proves the result. \endproof

%------------------------------------------------------------------------------
\begin{lemma}\label{TLfinal} For $n$ large enough, we have \hbox{$\w_n =\w^*_n$}. \end{lemma}
%------------------------------------------------------------------------------
\proof Let $\chi_n=\partial\w_n/\partial\theta$. Clearly $\int^\pi_{-\pi}\chi_n(t,\theta)\,d\theta=0$, and since $\w_n\in H^1(\mathcal C)$, then $\chi_n\in L^2(\mathcal C)$. Moreover, $\chi_n$ satisfies
\[\label{F3}-\Delta\chi_n+a^2_n\,\chi_n=(p_n-1)\,\big(\w_n(t,\theta)\big)^{p_n-2}\,\chi_n
\]
(in the sense of distributions), where
\begin{eqnarray*}\Big|(p_n-1)\,\big(\w_n(t_n,\theta)\big)^{p_n-2}\chi_n\Big|&\leq\!&(p_n-1)\,(\w_n(0,0))^{p_n-2}\Big(\frac{\w_n(t,\theta)}{\w_n(0,0)}\Big)^{p_n-2}|\chi_n|\\
&\leq\!&(p_n-1)\,(\w_n(0,0))^{p_n-2}\,|\chi_n|\in L^2(\mathcal C)\;.
\end{eqnarray*}
In other words, $-\Delta\chi_n+a^2_n\,\chi_n\in L^2(\mathcal C)$, and hence $\chi_n\in H^1(\mathcal C)$ satisfies:
$$
\|\nabla\chi_n\|^2_{L^2}+a^2_n\,\|\chi_n\|^2_{L^2}=(p_n-1)\,\int_\mathcal C\Big(\frac{\w_n(t,\theta)}{\w_n(0,0)}\Big)^{p_n-2}\chi^2_n\,dx\;.
$$
By Proposition~\ref{Proposition3.1}, we know that if $\psi\in H^1(\mathcal C)$ and $\int^\pi_{-\pi}\psi(t,\theta)\,d\theta=0$ a.e. $t\in\R$, then
\[
\|\nabla\psi\|^2_2-\beta_n\int_\mathcal C\frac{|\psi(t,\theta)|^2}{\big(\cosh((\alpha_n+1)\,t)\big)^2}\;dt\,d\theta\geq\Big[1-\(\frac{a_n\,p_n}2\)^2\Big]\|\psi\|^2_{L^2(\mathcal C)}
\]
with $\beta_n=a_n^2\,p_n\,(p_n-1)/2$. Passing to the limit as $n\to+\infty$, we get
\[
\|\nabla\psi\|^2_2-2\,(\alpha+1)^2\int_\mathcal C\frac{|\psi(t,\theta)|^2}{\big(\cosh((\alpha+1)\,t)\big)^2}\;dt\,d\theta\geq\Big[1-(\alpha+1)^2\Big]\|\psi\|^2_{L^2(\mathcal C)}\;.
\]
Consequently, for $\psi=\chi_n$, we obtain
\[
\begin{array}{rcl}
\hspace*{12pt}0&\hspace*{-6pt}=&\|\nabla\chi_n\|^2_2+a_n^2\,\|\chi_n\|^2_{L^2}-(p_n-1)\,\int_\mathcal C\big(\w_n(t,\theta)\big)^{p_n-2}\chi^2_n\,dx \vspace*{8pt}\\
&\hspace*{-6pt}=&\hspace*{-6pt}{\|\nabla\chi_n\|^2_{L^2}-2\,(\alpha+1)^2\int_\mathcal C\frac{\chi^2_n}{\(\cosh((\alpha+1)\,t)\)^2}\,dx+a^2_n\,\|\chi_n\|^2_{L^2(\mathcal C)}}\\
&&\hspace*{-6pt}+(p_n-1)\,\big(\w_n(0,0)\big)^{p_n-2}\displaystyle\int_\mathcal C\Big[{\textstyle\frac{1}{\(\cosh((\alpha+1)\,t)\)^2}-\Big(\frac{\w_n(t,\theta)}{\w_n(0,0)}\Big)^{p_n-2}}\Big]\chi^2_n\,dx\\
&&\hspace*{-6pt}+\big[2\,(\alpha+1)^2-(p_n-1)\,(\w_n(0,0))^{p_n-2}\big]\displaystyle\int_\mathcal C\frac{\chi^2_n}{\(\cosh((\alpha+1)\,t)\)^2}\,dx\\
&\hspace*{-6pt}\geq&\hspace*{-6pt}\big[1+a_n^2-(\alpha+1)^2-(p_n-1)\,\big(\w_n(0,0)\big)^{p_n-2}\,r_n\big]\,\|\chi_n\|^2_{L^2(\mathcal C)}\\
&&\hspace*{-6pt}+\big[2\,(\alpha+1)^2-(p_n-1)\,(\w_n(0,0))^{p_n-2}\big]\displaystyle\int_\mathcal C\frac{\chi^2_n}{\(\cosh((\alpha+1)\,t)\)^2}\,dx
\end{array}
\]
with $r_n=\sup_\mathcal C\big|\big(\w_n(t,\theta)/\w_n(0,0)\big)^{p_n-2}\kern-12pt-e^V\big|$. Recall that by Lemma~\ref{Lem:Chou-Wan},
$$\lim_{n\to+\infty}(p_n-1)(\w_n(0,0))^{p_n-2}=\mu=2\,(\alpha+1)^2\;,
$$
and by Lemma~\ref{Proposition1}, $\lim_{n\to+\infty}r_n=0$. Since $a_n\to0$ as $n\to+\infty$ and $(1+\alpha)^2<1$, we readily get a contradiction for large $n$, unless $\chi_n\equiv 0$. This means that $\w_n$ is independent of the variable~$\theta$, and so $\w_n=\w^*_n$. \endproof

%%%%%%%%%%%%%%%%%%%%%%%%%%%%%%%%%%%%%%%%%%%%%%%%%%%%%%%%%%%%%%%%%%%%%%%%%%%%%%%
%%%%%%%%%%%%%%%%%%%%%%%%%%%%%%%%%%%%%%%%%%%%%%%%%%%%%%%%%%%%%%%%%%%%%%%%%%%%%%%
\section{Concluding remarks}\label{sect5}

It is interesting to note that, via the Emden-Fowler transformation \eqref{3.1}, for any $\alpha>-1$, inequality \eqref{MTalpha} can be stated on the space
\[
\mathfrak E_\alpha=\Big\{w=w(t,\theta)\in L^1(\mathcal C,d\nu_\alpha)\;:\;|\nabla w|\in L^2 (\mathcal C,dx)\Big\}
\]
where
\[
d\nu_\alpha:=\frac{\alpha+1}2\,\frac{dt}{\big[\cosh\big((\alpha+1)\,t\big)\big]^2}\;.
\]

%------------------------------------------------------------------------------
\begin{proposition}\label{Prop:MTalphaCylinder} If $\alpha>-1$, then
\[\label{MTalphaCylinder}
\int_{\mathcal C} e^{w-\int_{\mathcal C}w\,d\nu_\alpha}\;d\nu_\alpha\,\leq\,e^{ \frac{1}{16\,\pi\,(\alpha+1)}\,\left(\|\nabla w\|^2_{L^2(\mathcal C)} +\alpha\,(\alpha+2)\,\|\,\partial_\theta w\,\|^2_{L^2(\mathcal C)}\right)}\quad\forall\;w\in \mathfrak E_\alpha\;.
\]
\end{proposition}

\bigskip

\medskip

%------------------------------------------------------------------------------
As in Section~\ref{sect2.1}, when $\alpha\leq 0$, there holds
$$
\int_{\mathcal C} e^{w-\int_{\mathcal C}w\,d\nu_\alpha}\;d\nu_\alpha\,\leq\,e^{ \frac{1}{16\,\pi\,(\alpha+1)}\,\|\nabla w\|^2_{L^2(\mathcal C)}}\quad\forall\;w\in \mathfrak E_\alpha\;,
$$
with extremals obtained from \eqref{ExtremalCylinder} up to translations, scalings and addition of constants. 

However, when $\alpha>0$, while the latter inequality is always valid for functions depending only on the variable $t\in\R$, in general it fails to hold in $\mathfrak E_\alpha$.

\medskip The above inequality is one of the three equivalent versions of the weighted Moser-Trudinger inequalities that we prove in this paper: on the sphere $S^2$, on the euclidean space $\R^2$ and on the cylinder $\mathcal C$. The symmetry breaking phenomenon is easily understood in this case, as clearly, the corresponding extremals are symmetric if and only if $\alpha\in(-1,0]$.

On the contrary, the symmetry breaking phenomenon in Caffarelli-Kohn-Nirenberg inequality is a more subtle issue, since it is less evident how the weights conspire against symmetry. Our key observation is that weighted Moser-Trudinger inequalities appear as limits of Caffarelli-Kohn-Nirenberg inequalities in an appropriate blow-up limit. In this asymptotics, the case $b<h(a)$ yields to $\alpha>0$, while the case $b>h(a)$ leads to $\alpha\in(-1,0)$.

%%%%%%%%%%%%%%%%%%%%%%%%%%%%%%%%%%%%%%%%%%%%%%%%%%%%%%%%%%%%%%%%%%%%%%%%%%%%%%
%%%%%%%%%%%%%%%%%%%%%%%%%%%%%%%%%%%%%%%%%%%%%%%%%%%%%%%%%%%%%%%%%%%%%%%%%%%%%%
\section*{Appendix. The dilated stereographic projection}

We use spherical coordinates $(\phi,\theta)\in[-\frac\pi 2,\frac\pi 2]\times[0,2\pi)$ on $S^2\subset\R^3$ and radial coordinates $(r,\theta)\in[0,\infty)\times[0,2\pi)$ on $\R^2$. By definition of the dilated stereographic projection, we have
\[
\cos\phi=\frac{2\,r^{\alpha+1}}{1+r^{2(\alpha+1)}}\quad\mbox{and}\quad\sin\phi=\frac{r^{2(\alpha+1)}-1}{1+r^{2(\alpha+1)}}\;,
\]
from which we deduce
\[
\cos\phi\,\frac{d\phi}{dr}=\frac{4\,(\alpha+1)\,r^{2\alpha+1}}{(1+r^{2(\alpha+1)})^2}\;.
\]
The normalized measure of the sphere $S^2$ is given by
\[
d\sigma=\frac 12\,\cos\phi\,\frac{d\theta}{2\pi}
\]
and a simple change of variables shows that, if $u(\phi,\theta)=v(r,\theta)$, then
\[
\int_{S^2}f(u)\;d\sigma=\int_{\R^2}f(v)\,\frac{\cos\phi}2\,\frac{d\phi}{dr}\;dr\,\frac{d\theta}{2\pi}=\int_{\R^2}f(v)\,d\mu_\alpha
\]
where $d\mu_\alpha=\frac{\alpha+1}\pi\,\frac{r^{2\alpha}}{(1+r^{2(\alpha+1)})^2}\;r\,dr\,d\theta$. Using spherical and radial coordinates respectively on $S^2$ and $\R^2$, the expressions of the gradients are given respectively as follows
\[
|\nabla u|^2=|\partial_\phi u|^2+\frac 1{\cos^2\phi}\,|\partial_\theta u|^2\quad\mbox{and}\quad |\nabla v|^2=|\partial_r v|^2+\frac 1{r^2}\,|\partial_\theta r|^2\;.
\]
Knowing that $\partial_\phi u=\partial_rv\,\left(\frac{d\phi}{dr}\right)^{-1}$, we get
\[
\int_{S^2}\!|\partial_\phi u|^2\,d\sigma=\!\int_{\R^2}\!|\partial_rv|^2\,\frac{\cos\phi}2\left(\frac{d\phi}{dr}\right)^{\!-1}\!\!dr\,\frac{d\theta}{2\pi}=\frac 1{4\pi\,(\alpha+1)}\int_{\R^2}\!|\partial_rv|^2\,r\,dr\,d\theta\;.
\]
While using that $\partial_\theta u=\partial_\theta v$, we get
\[
\int_{S^2}\frac 1{\cos^2\phi}\,|\partial_\theta u|^2\,d\sigma=\!\int_{\R^2}\!|\partial_\theta v|^2\,\frac 1{2\cos\phi}\,\frac{d\phi}{dr}\;dr\,\frac{d\theta}{2\pi}=\frac{\alpha+1}{4\pi}\!\int_{\R^2}\!\frac{|\partial_\theta v|^2}{r^2}\,r\,dr\,d\theta\;.
\]
Thus,  observing that $(\alpha+1)^2-1=\alpha\,(\alpha+2)$, we conclude
\[
\int_{S^2}|\nabla u|^2\,d\sigma=\frac 1{4\pi\,(\alpha+1)}\left[\int_{\R^2}|\nabla v|^2\,dx + \alpha\,(\alpha+2)\int_{\R^2}\frac{|\partial_\theta v|^2}{r^2}\,r\,dr\,d\theta\right]\;.
\]

%%%%%%%%%%%%%%%%%%%%%%%%%%%%%%%%%%%%%%%%%%%%%%%%%%%%%%%%%%%%%%%%%%%%%%%%%%%%%%
%%%%%%%%%%%%%%%%%%%%%%%%%%%%%%%%%%%%%%%%%%%%%%%%%%%%%%%%%%%%%%%%%%%%%%%%%%%%%%

\bigskip\noindent{\small{\bf Acknowlegments.} This work has been partially supported by European Programs HPRN-CT \#~2002-00277 \& 00282, by the projects ACCQUAREL and IFO of the French National Research Agency (ANR) and by M.U.R.S.T. project: Variational Methods and Non Linear Differential Equation, Italy. The third author wishes also to express her gratitude to Ceremade for the warm and kind hospitality during her visits.}

\medskip\noindent{\small \copyright\,2007 by the authors. This paper may be reproduced, in its entirety, for non-commercial purposes.}

%%%%%%%%%%%%%%%%%%%%%%%%%%%%%%%%%%%%%%%%%%%%%%%%%%%%%%%%%%%%%%%%%%%%%%%%%%%%%%%
%\nocite*\bibliographystyle{siam}\bibliography{References.bib}

\begin{thebibliography}{10}

\bibitem{Beckner-93}
{\sc W.~Beckner}, {\em Sharp {S}obolev inequalities on the sphere and the
 {M}oser-{T}rudinger inequality}, Ann. of Math. (2), 138 (1993), pp.~213--242.

\bibitem{Brezis-Merle-91}
{\sc H.~Brezis and F.~Merle}, {\em Uniform estimates and blow-up behavior for
 solutions of {$-\Delta u=V(x)e\sp u$} in two dimensions}, Comm. Partial
 Differential Equations, 16 (1991), pp.~1223--1253.

\bibitem{Caffarelli-Kohn-Nirenberg-84}
{\sc L.~Caffarelli, R.~Kohn, and L.~Nirenberg}, {\em First order interpolation
 inequalities with weights}, Compositio Math., 53 (1984), pp.~259--275.

\bibitem{Catrina-Wang-01}
{\sc F.~Catrina and Z.-Q. Wang}, {\em On the {C}affarelli-{K}ohn-{N}irenberg
 inequalities: sharp constants, existence (and nonexistence), and symmetry of
 extremal functions}, Comm. Pure Appl. Math., 54 (2001), pp.~229--258.

\bibitem{MR1121147}
{\sc W.~X. Chen and C.~Li}, {\em Classification of solutions of some nonlinear
 elliptic equations}, Duke Math. J., 63 (1991), pp.~615--622.

\bibitem{Chou-Chu-93}
{\sc K.~S. Chou and C.~W. Chu}, {\em On the best constant for a weighted
 {S}obolev-{H}ardy inequality}, J. London Math. Soc. (2), 48 (1993),
 pp.~137--151.

\bibitem{Chou-Wan-94}
{\sc K.~S. Chou and T.~Y.-H. Wan}, {\em Asymptotic radial symmetry for
 solutions of {$\Delta u+e\sp u=0$} in a punctured disc}, Pacific J. Math.,
 163 (1994), pp.~269--276.

\bibitem{MR1372247}
{\sc K.~S. Chou and T.~Y.~H. Wan}, {\em Correction to: ``{A}symptotic radial
 symmetry for solutions of {$\Delta u+e\sp u=0$} in a punctured disc''
 [{P}acific {J}. {M}ath.\ {\bf 163} (1994), no.\ 2, 269--276]}, Pacific J.
 Math., 171 (1995), pp.~589--590.

\bibitem{Felli-Schneider-03}
{\sc V.~Felli and M.~Schneider}, {\em Perturbation results of critical elliptic
 equations of {C}affarelli-{K}ohn-{N}irenberg type}, J. Differential
 Equations, 191 (2003), pp.~121--142.

\bibitem{Landau-Lifschitz-67}
{\sc L.~Landau and E.~Lifschitz}, {\em Physique th\'eorique. Tome III:
 M\'ecanique quantique. Th\'eorie non relativiste. (French)}, Deuxi\`eme
 \'edition. Translated from russian by E. Gloukhian. \'Editions Mir, Moscow,
 1967.

\bibitem{MR2053993}
{\sc C.-S. Lin and Z.-Q. Wang}, {\em Erratum to: ``{S}ymmetry of extremal
 functions for the {C}affarelli-{K}ohn-{N}irenberg inequalities'' [{P}roc.\
 {A}mer.\ {M}ath.\ {S}oc.\ {\bf 132} (2004), no.\ 6, 1685--1691]}, Proc. Amer.
 Math. Soc., 132 (2004), p.~2183 (electronic).

\bibitem{Lin-Wang-04}
\leavevmode\vrule height 2pt depth -1.6pt width 23pt, {\em Symmetry of extremal
 functions for the {C}affarelli-{K}ohn-{N}irenberg inequalities}, Proc. Amer.
 Math. Soc., 132 (2004), pp.~1685--1691 (electronic).

\bibitem{MR2001882}
{\sc D.~Smets and M.~Willem}, {\em Partial symmetry and asymptotic behavior for
 some elliptic variational problems}, Calc. Var. Partial Differential
 Equations, 18 (2003), pp.~57--75.

\bibitem{Tarantello-07}
{\sc G.~Tarantello}, {\em Selfdual gauge field vortices: an analytical
 approach}, PNLDE 72, Birkh{\"a}user ed. Boston MA, USA, 2007.

\end{thebibliography}

\end{document}